%% file: main.tex
\documentclass[11pt,reqno,oneside]{amsart}
\calclayout
\usepackage{hyperref}
\usepackage[foot]{amsaddr}
\usepackage[T1]{fontenc}
\usepackage[english]{babel}
\addto\extrasenglish{
  
}
\usepackage{amssymb}
\usepackage{amsmath}
\usepackage{amsthm}
\usepackage{amsfonts}
\usepackage[shortlabels]{enumitem}
\usepackage{graphicx}
\usepackage{mathabx}
\usepackage{eurosym}
\usepackage{dsfont}
\usepackage{csquotes}
\usepackage{xparse}
\usepackage{accents}
\usepackage{color}
\usepackage{xcolor}
\usepackage{blindtext}
\usepackage{thmtools}
\usepackage{mathrsfs}
\usepackage{stmaryrd}
\usepackage{thm-restate}
\usepackage[
backend=biber,
style=alphabetic,
]{biblatex}

\usepackage{xurl}
\hypersetup{breaklinks=true}
\hypersetup{
hypertexnames=false,
colorlinks = true,
linkcolor = {purple},
urlcolor = {blue},
citecolor = {brown}
}
\newtheorem{lemma}{Lemma}[section]

\newtheorem{definition}[lemma]{Definition}
\newtheorem{theorem}[lemma]{Theorem}

\newtheorem{proposition}[lemma]{Proposition}
\newtheorem{ass}{A\!\!}
\newtheorem{manualassinner}{A\!\!}
\newenvironment{manualass}[1]{%
  \renewcommand\themanualassinner{#1}%
  \manualassinner
}{\endmanualassinner}
\newtheorem{example}[lemma]{Example}
\newtheorem{remark}[lemma]{Remark}
\makeatletter \@addtoreset{equation}{section} \makeatother

\newcommand{\RR}{{\mathbb R}}
\newcommand{\NN}{{\mathbb N}}
\newcommand{\ZZ}{{\mathbb Z}}
\newcommand{\lun}{{\mathcal{M}(\XX)}}
\newcommand{\linf}{{\mathcal{B}(\XX)}}

\newcommand{\EE}{\mathbb{E}}
\newcommand{\E}{\mathcal{E}}

\newcommand{\XX}{\mathbb{X}}
\newcommand{\X}{\mathcal{X}}
\newcommand{\Q}{\mathcal{Q}}

\newcommand*{\A}{\mathcal A}

\newcommand{\M}{\mathcal{M}}
\renewcommand{\P}{\mathcal{P}}
\newcommand{\B}{\mathcal{B}}
\newcommand{\MX}{\mathcal{K}^+}

\newcommand{\bc}{(1-c)}
\newcommand{\diam}{diam}

\newcommand{\mat}[2]{#1(#2)}

\newcommand{\stardot}{\overset{\cdot}{\star}}
\DeclareMathOperator{\Id}{Id}
\let\diam\relax
\DeclareMathOperator{\diam}{\mathbf{diam}}
\newcommand{\PP}{{\mathbb P}}

\newcommand\arr[2]{\left\{ #1 \right\}_{#2}}
\def\beq{\begin{equation}}
\def\eeq{\end{equation}}
\let\origautoref\autoref
\def\autoref#1{\textbf{\origautoref{#1}}}

\title{Ergodic behavior of products of random positive operators}
\author{Maxime Ligonnière}
\address{Institut Denis Poisson UMR 7013, Université de Tours, Université d’Orléans, CNRS France}
\email{maxime.ligonniere@lmpt.univ-tours.fr}
\address{Ecole Polytechnique, Centre de mathématiques appliquées (CMAP), 91128 Palaiseau, France}

\begin{document}

\begin{abstract}
This article is devoted to the study of products of random operators of the form $M_{0,n}=M_0\cdots M_{n-1}$, where $(M_{n})_{n\in\NN}$ is an ergodic sequence of positive operators on the space of signed measures on a space $\XX$. Under suitable conditions, in particular, a Doeblin-type minoration suited for non conservative operators, we obtain asymptotic results of the form 
\[  \mu M_{0,n} \simeq \mu(\tilde{h}) r_n \pi_n,\]
where $\tilde{h}$ is a random bounded function, $(r_n)_{n\geq 0}$ is a random non negative sequence and $\pi_n$ is a random probability measure on $\XX$. Moreover, $\tilde{h}$, $(r_n)$ and $\pi_n$ do not depend on the choice of the measure $\mu$. We prove additionally that $n^{-1} \log (r_n)$ converges almost surely to the Lyapunov exponent $\lambda$ of the process $(M_{0,n})_{n\geq 0}$ and that the sequence of random probability measures $(\pi_n)$ converges weakly towards a random probability measure. These results are analogous to previous estimates from \cite{hennion_limit_1997} in the case of $d\times d$ matrices, that were obtained with different techniques, based on a projective contraction in Hilbert distance. 
In the case where the sequence $(M_n)$ is i.i.d, we additionally exhibit an expression of the Lyapunov exponent $\lambda$ as an integral with respect to the weak limit of the sequence of random probability measures $(\pi_n)$ and exhibit an oscillation behavior of $r_n$ when $\lambda=0$.
We provide a detailed comparison of our assumptions with the ones of \cite{hennion_limit_1997} and present some example of applications of our results, in particular in the field of population dynamics.
\end{abstract}
\maketitle
\setcounter{tocdepth}{1}
\tableofcontents

\input{1_Intro.tex}
\input{2_results.tex}

\input{3_proofs.tex}
\input{4_hennion.tex}
\input{5_branching.tex}

\section{Acknowledgements}
I have received the support of the Chair "Modélisation Mathématique et Biodiversité" of VEOLIA-Ecole Polytechnique-MnHn-FX, and of the ANR project NOLO (ANR 20-CE40-0015), funded by the French ministry of research. I would like to warmly thank my PhD-supervisors Vincent Bansaye and Marc Peigné for their continous guidance and support as well as their numerous feedbacks on this manuscript, which greatly helped improving it. I would also like to thank Gerold Alsmeyer for a fruitful discussion on the topic of Markov Random Walks.

\printbibliography

\end{document}

%% file: 1_Intro.tex
\section{Introduction}
\subsection{General introduction}
The study of products of random linear operators can be traced back to the seminal article of \cite{furstenberg_products_1960}, studying products of the form 
\[M_{0,n}=M_0\dots M_{n-1},\]
where $(M_{n})_{n\geq 0}$ is a stationary and metrically transitive sequence of $p\times p$ real or complex random matrices. Under a mild irreducibility assumption, the authors exhibit a strong law of large numbers on the norm of $M_{0,n}$, in the form 
\[\lim_{n\rightarrow \infty} \frac{1}{n} \log \Vert \mu M_{0,n} \Vert=\lambda, \]
where $\lambda$ is a deterministic number called Lyapunov exponent of the sequence $(M_{n})_{n\geq 0}$, defined as \[ \lambda = \underset{n\rightarrow \infty}{\lim}\frac{1}{n} \EE \left[  \log \Vert M_{0,n} \Vert \right]= \inf_{n\geq 1} \frac{1}{n} \EE \left[ \log \Vert M_{0,n} \Vert \right],\] where the norm $\Vert \cdot \Vert$ can be chosen to be any submultiplicative norm.  Under additional positivity and boundedness assumptions on the entries of the matrices $(M_n)$, \cite{furstenberg_products_1960} also prove a strong law of large numbers for the entries $M_{0,n}(i,j)$ of the products : for any $i,j\in\{1,\dots,d\}$, almost surely, \[ \underset{n\rightarrow \infty}\lim \frac{1}{n} \log \mat{M_{0,n}}{i,j}=\lambda.\]
These estimates on the behavior of the entries of $M_{0,n}$ were then extended to the case of products of invertible matrices, see e.g. \cite{guivarch_proprietes_2001} and \cite{bougerol_products_1985}. These works rely on a careful study of the action of invertible matrices on the projective space $\P(\RR^d)$.
\\ To strengthen the results from \cite{furstenberg_products_1960} on products of matrices with non negative entries and relax their assumptions, \cite{hennion_limit_1997} studied the action of $\M_d(\RR_+)$ on the projective space $\P(\RR_+^d)$, endowed with the so called (pseudo)-Hilbert distance $d_H$ previously defined in \cite{busemann_projective_1953} and \cite{birkhoff_extensions_1957}. This distance is particularly well adapted to this problem, since the contraction coefficient of the projective action of a matrix with respect to $d_H$ is explicit in terms of its entries, in particular, any matrix with nonnegative entries in $1$-contracting and any matrix with positive entries is strictly contracting. 
Under the assumption that almost surely, for $n$ large enough, $M_{0,n}$ has all positive entries, Hennion obtains the asymptotic decomposition 
\[ M_{0,n}(i,j)=\lambda_n R_n(i) L_n(j) + \underset{n\rightarrow \infty}o(\lambda_n),\]
where $\lambda_n$ is the dominant eigenvalue of $M_{0,n}$ and $L_n,R_n$ are the associated left and right eigenvectors, with the normalizations $\Vert R_n \Vert =1$ and $\langle L_n, R_n \rangle=1$. Moreover $(R_n)_{n\geq 0}$ almost surely converges to a random vector $R$, $\left(L_n/\Vert L_n \Vert \right)_{n} $ converges in distribution, and $\left(n^{-1} \log \lambda_n\right)_{n\geq 1}$ almost surely converges to the Lyapunov exponent $\lambda$.
\\ Such results have important implications, in particular in the field of populations dynamics. Indeed, a population composed of $d$ types of individuals, evolving in a fluctuating environment, without interacting which each other, can be modelled by a linear model of the form \begin{equation} \label{eq: linear model} x_n= x_{n-1} M_{n-1},\end{equation}  where $x_n$ is a row vector of $\RR_+^d$ encoding the mass of individuals of each type at time $n$ and $M_{n-1}=\left(M_{n-1}(i,j)\right)_{1\leq i,j\leq d}$ is a random matrix encoding the rates at which individuals of each type $i$ create individuals of each type $j$ between times $n-1$ and $n$. In such a time-inhomogeneous population model, the understanding of the asymptotics of $x_n$ amounts to the understanding of the matrix product $M_{0,n}$.
\\Moreover, such products also appear in the study of multitype Galton-Watson processes in random environment (MGWRE), which were introduced in \cite{athreya_branching_1971}. They are a generalization of Galton-Watson processes to the case where the distribution of the (random) offspring of an individual depends on a notion of type and on a random environment that changes through time. When conditioning a MGWRE on the environment sequence, one obtains a so-called quenched population model, which satisifies \eqref{eq: linear model}, where $x_n$ is the expectation of the population conditionally on the environmental sequence. The value of the Lyapunov exponent $\lambda$ of the underlying matrix product separates three regimes of the MGWRE : subcritical ($\lambda <0$), critical ($\lambda=0$), supercritical ($\lambda>0$). These three regimes have different properties. In particular, when $\lambda\leq 0$, the MGWRE goes extinct with probability $1$, when $\lambda>0$, the MGWRE survives with positive probability. This separation between regimes was established in \cite{athreya_branching_1971} and \cite{kaplan_results_1974}, using results from \cite{furstenberg_products_1960}. More recent advances in the study of random matrix products - in particular Hennion's article- were key to the last developments of the theory of MGWRE in random environments, see e.g. \cite{cam_pham_conditioned_2018, le_page_survival_2018, grama_kestenstigum_2023}.
\\ Products of random infinite dimensional operators have also been the subject of some investigation. In the case where the involved operators are the transition matrices of some Markov chain, in other words if they are conservative (in the sense that $M\mathds{1}=\mathds{1}$), ergodicity results are obtained by \cite{cogburn_ergodic_1984}, completed with some more precise results in \cite{orey_markov_1991}. Various limit theorems (law of large numbers, existence of Lyapunov exponents, central limit theorem, local central limit theorem, Large Deviation principle) have also been obtained on products of ergodic sequences of infinite dimensional operators using spectral techniques by \cite{dragicevic_spectral_2018}. The stability of the Lyapunov exponents under perturbation is studied in \cite{atnip_perturbation_2022,froyland_hilbert_2019}. Ergodicity of products of random operators was also obtained in \cite{kifer_perron-frobenius_1996} in the case where they act on some compact space $\XX$.
\\ In this paper, we would like to obtain similar ergodicity results for products of infinite dimensional positive operators, thus extending the results of \cite{kifer_perron-frobenius_1996} without any topological assumption on $\XX$, in particular without compacity. We have in particular in mind applications to population models with an infinite number of types. We first consider such a set $\XX$, typically infinite, endowed with a $\sigma$-algebra $\X$, and build a set $\MX$ of positive linear operators acting both on the space of signed measures $\mathcal{M}(\XX)$ on the left and the space of measurable bounded functions $\linf$ on the right. Then, we let $(M_n)$ be a stationary, ergodic sequence of elements of $\MX$ and define the products $M_{0,n}=M_0\cdots M_{n-1}$. The approach of \cite{hennion_limit_1997} can be extended to this infinite dimensional setup. Indeed, it is possible to define the Hilbert distance $d_H$ on the projective positive cone of an infinite dimensional vector space and to obtain a nice characterisation of the operators that are (strictly) contracting with respect to $d_H$. We refer the reader to \cite{ligonniere_contraction_2023} for a proof of these facts. However, as we explain in Section \ref{sec: Hilbert}, this characterisation leads to stronger positivity assumptions in an infinite dimensional context than it did in the finite dimensional one. For example, such an extension of Hennion's approach would not be able to deal with products of infinite Leslie matrices that we present in Section \ref{sec: branching models}. 
\\ For this reason, we use a different contraction method to obtain a projective contraction. This method aims at extending the Doeblin contraction techniques for Markov operators to a product of non conservative operators (that is operators $M$ such that $M\mathds{1}\neq \mathds{1}$ in general). To do so, we consider the auxiliary family of Markov operators $P_{k,n}^N$, defined for each $k\leq n \leq N$ as\[\delta_x P_{k,n}^N f = \frac{ \delta_x M_{k,n} (f m_{n,N})}{m_{k,N}(x)},\]
where, for $x\in\XX$, \[m_{k,n}(x)=\delta_x M_{k,n}\mathds{1}.\]
These Markov operators are related to the projective action $\mu\cdot M_{k,n}=\frac{\mu M_{k,n}}{\Vert \mu M_{k,n} \Vert}$ of $M_{k,n}$ on measures by :
\[ \delta_x \cdot M_{k,n} = \delta_x P_{k,n}^n = \delta_x P_{k,k+1}^n\cdots P_{n-1,n}^n,\]
for any $x\in\XX$.
We provide sufficient conditions for the Markov operators $(P_{k,k+1}^n)_{k< n}$ to satisfy Doeblin minorations of the form $\delta_x P(f) \geq c \nu(f)$, which guarantees that they are contracting in total variation. This allows to obtain a notion of projective contraction of the $M_{k,n}$ on the set of signed measures. Such auxiliary operators were already introduced and already applied to study both homogeneous semi-groups of operators, e.g. in  \cite{del_moral_stability_2002, champagnat_exponential_2016} and inhomogeneous ones, e.g. in \cite{bansaye_ergodic_2020}. We consider here a random sequence of operators $(M_n)$, i.e a discrete time, random, time-inhomogeneous semi group, and assume this sequence is stationary and ergodic. The stationary and ergodic framework allows us to provide more explicit assumptions, and we obtain a finer asymptotic analysis of the sequence $(M_{0,n})_{n\geq0}$. Namely, we prove that the following almost sure approximation in total variation holds for $n$ large enough
\begin{equation} \label{eq: presentation resultat} \left \Vert \mu M_{0,n} - \mu(h) r_n \pi_n \right\Vert_{TV}\leq \delta^n\Vert \mu M_{0,n}\Vert_{TV},  \end{equation}
where $\delta<1$, ${h}$ is a random bounded function, $r_n$ is a positive random number  and $(\pi_n)$ is a sequence of random probability measures on $\XX$, which are all independent of the measure $\mu$. We prove additionally that $(n^{-1} \log (r_n))$ converges almost surely to the Lyapunov exponent of the process $M_{0,n}$, and that the sequence $(\pi_n)$ of random probabilities converges in distribution with respect to the total variation topology towards a random probability measure $\Lambda$ on $\mathcal{M}_1(\XX)$.
Additionally, we show in Theorem \ref{thme:lyap} that when the sequence $(M_n)$ is i.i.d, the probability distribution $\Lambda$ and the Lyapunov exponent $\lambda$ are related as follows :
\[ \lambda=\int \log \Vert \mu M \Vert d\Lambda(\mu) d\P(M),\]
where $\P$ refers to the law of the operators $(M_n)$, thus extending a result stated in \cite{bougerol_products_1985} in finite dimension.
Finally, still under the assumption that $(M_n)$ is i.i.d, we show in Theorem \ref{thme:osc} that, when $\lambda=0$, it holds almost surely, for any $\mu\in\mathcal{M}_+(\XX)-\{0\}$,
\[\tag{OSC} \underset{n\rightarrow\infty}\limsup \log \Vert \mu M_{0,n} \Vert = -\underset{n\rightarrow\infty}\liminf \Vert \mu M_{0,n}\Vert = +\infty,\]
\stepcounter{equation}
except in a situation knowed as Null-Homology.
\\These results should allow to extend many known results on MGWRE with a finite type set to a class of MGWRE with an infinite type set $\XX$. In particular, our results imply that when the Lyapunov exponent $\lambda$ is nonpositive, outside of Null-Homology, the quenched population size $ \mu M_{0,n} \mathds{1}$ satisfies $\underset{n\rightarrow \infty}\liminf \mu M_{0,n}\mathds{1}=0$ almost surely. By a classical first moment argument, this is a sufficient condition for the almost sure extinction of the population. The survival of the population when $\lambda>0$ is a more delicate problem and will be the object of a forecoming article. 
\subsection{Framework and notations}
Let $(\XX,\X)$ be a measurable set of arbitrary cardinality, such that for all $x\in \XX$, the singleton $\{ x \} \in \X$. We denote by $\mathcal{B}(\XX)$ the Banach space of bounded measurable functions on $\XX$, endowed with the supremum norm, and $\mathcal{B}_+(\XX)$ the cone of nonnegative functions of $\mathcal{B}(\XX)$. The vector space of signed measures, noted $\M(\XX)$, and the cone of nonnegative elements of $\M(\XX)$, noted $\mathcal{M}_+(\XX)$, are endowed with the total variation norm $\Vert \cdot \Vert_{TV}$. Note that $\Vert \mu \Vert_{TV}=\mu(\XX)$ for any $\mu \in \mathcal{M}_+(\XX)$. Let $\mathcal{M}_1(\XX)$ be the set of probability measures on $(\XX,\X)$. For any measurable set $A\in \X$, let $\mathds{1}_A\in\mathcal{B}(\XX)$ be the indicator function of $A$. For short, we note $\mathds{1}=\mathds{1}_\XX\in \mathcal{B}(\XX)$ the unit function on $\XX$. We also note $\delta_x\in \mathcal{M}_1(\XX)$ the Dirac measure at $x$.

Let $\MX$ be the set of maps $Q$ of the form:
\[ Q : \XX \times \X \longrightarrow \RR^+,\]
such that, for any $x\in\XX$, for any $A \in \X$, the map $x\mapsto Q(x,A)$ is measurable, the map $A\mapsto Q(x,A)$ is a positive and finite measure on $(\XX,\X)$ and $\vvvert Q \vvvert:= \sup_{x\in\XX} Q(x,\XX)<\infty$.

Such a map $Q\in\MX$ naturally operates on $\mathcal{B}(\XX)$ by setting, for any $f\in \mathcal{B}(\XX)$, and any $x\in\XX$,
\[Qf(x)= \int_{\XX} f(y) Q(x,dy).\]
Note that $\vert Qf(x)\vert \leq \Vert f \Vert_\infty \vvvert Q \vvvert$, thus $Qf\in \mathcal{B}_+(\XX)$ as soon as $f\in\mathcal{B}_+(\XX)$ and $Q$ acts as a bounded positive operator with norm $\vvvert Q \vvvert$ on $\mathcal{B}(\XX)$.
Moreover, for any positive measure $\mu\in \mathcal{M}_+(\XX)$, and any $Q\in\MX$, the positive measure $\mu Q$ on $\XX$ is well defined by setting, for any nonnegative function $f$,
\[ \mu Q(f) = \mu(Qf)= \int_{\XX} Qf(x) \mu(dx).\]
Note that $\mu Q$ has indeed finite mass $\mu Q(\mathds{1})=\mu(Q\mathds{1})\leq \vvvert Q \vvvert \mu(\mathds{1})<\infty$ since $\mu$ is assumed to be a finite measure. 
This action can therefore naturally be extended to the set of signed measures $\M(\XX)$, where $Q$ acts as a bounded linear operator with norm $\vvvert Q \vvvert$.

Thus, the elements of $\MX$ operate as positive linear operators both on the sets of bounded measurable functions and on the set of signed measures on $\XX$, with a duality relation between these two actions. 
Moreover, it is also possible to define a projective action $\cdot$ of $\MX$ onto the projective space associated with $\mathcal{M}_+(\XX)$, ie the set of probability measures $\mathcal{M}_1(\XX)$, by setting, for any $\mu\in \mathcal{M}_+(\XX)$ and any $M\in \MX$ such that $\mu M \neq 0$,
\[ \mu \cdot M = \frac{\mu M}{\Vert \mu M \Vert_{TV}}\in \mathcal{M}_1(\XX).\]
Finally, the set $\MX$ is naturally endowed with an associative, non commutative product, defined by : for any $Q_1,Q_2\in \MX$, any $x\in\XX$ and any $A\in \X$,
\[Q_1Q_2(x,A)= \int_{y} Q_1(x,dy)Q_2(y,A).\]
This product is compatible with the left and right actions defined above, in other words, for any $Q_1,Q_2\in \MX$, any $\mu \in \mathcal{M}_+(\XX)$, and $f\in \mathcal{B}_+(\XX)$,
\[\mu(Q_1Q_2)=(\mu Q_1) Q_2 \text{ and } Q_1 Q_2 (f)=Q_1(Q_2 f),\]
and whenever $\mu\in \mathcal{M}_+(\XX)$ and $\mu Q_1Q_2 \neq 0$, 
\[\mu \cdot (Q_1Q_2)=(\mu\cdot Q_1)\cdot Q_2.\]
The operator norm $\vvvert \cdot \vvvert $ satisfies the submultiplicativity relation \[ \vvvert Q_1 Q_2 \vvvert \leq \vvvert Q_1 \vvvert \vvvert Q_2 \vvvert.\]

\begin{remark}\label{remark : matrix}
In the case of a finite or countable set $\XX$, any measure on $\XX$ is atomic thus an operator $Q\in \MX$ on $\XX$ corresponds to a matrix indexed by $\XX$ with nonnegative entries. The product of operators of $\MX$ corresponds to the matrix product and the respective left and right actions of $\MX$ on signed measures and bounded functions correspond to the product of matrices respectively with the vectors of $\ell^1(\XX)$ (seen as row vectors) and of $\ell^\infty(\XX)$ (seen as column vectors).
\end{remark}
We consider a dynamical system $(\Omega, \mathcal{A}, \PP, \theta)$, where $(\Omega, \mathcal{A}, \PP)$ is a probability set and $\theta:\Omega\longrightarrow \Omega$ is a measurable transformation, which preserves the probability $\PP$, i.e  $\PP\circ \theta^{-1}=\PP$. 
\\Let $M:\Omega \longrightarrow \MX$ be a measurable map. We denote as $\NN_0$ the set of nonnegative integers and note, for each $n\in\NN_0$, 
\[M_n=M\circ \theta^n.\]
Note that the sequence $(M_n)_{n\in\NN_0}$ is stationary.
For each $k < n$, $\omega\in\Omega$, let us define the random product
\[M_{k,n}(\omega)= M_k(\omega)\cdots M_{n-1}(\omega)= (M\circ\theta^k(\omega))\cdots  \left(M\circ \theta^{n-1}(\omega)\right)\in \MX.\]
with the convention  $M_{k,k}(\omega)=\Id\in \MX$. Notice that \(M_{k,k+n}(\omega)=M_{0,n}\circ\theta^k (\omega).\) The operators satisfy the following semi group property : for any $k\leq n\leq N$, any $\mu \in \mathcal{M}(\XX),$ any  $f \in \mathcal{B}(\XX)$
\begin{equation} \label{eq: semigroup}
 \mu M_{k,N}(\omega) f = \mu M_{k,n}(\omega) M_{n,N}(\omega) f.
\end{equation}
Moreover, for any $x\in\XX$, $k\leq n$, $\omega\in\Omega$, we set \[m_{k,n}(x,\omega)= \delta_x M_{k,n}(\omega) \mathds{1}=\Vert\delta_x M_{k,n}(\omega) \Vert_{TV}.\]
Notice in particular that for any positive measure $\mu$,
\[ \Vert \mu M_{k,N} \Vert=\mu(m_{k,N})=\mu M_{k,N} \mathds{1}= \mu M_{k,n} m_{n,N}.\]
Let us point out additionally that \(\vvvert M_{k,n}(\omega)\vvvert=\sup_{x\in\XX} m_{k,n}(x,\omega)=\Vert m_{k,n}(\cdot,\omega) \Vert_\infty,\) and that for any $k\leq n \leq N$, \(\vvvert M_{k,N}\vvvert \leq \vvvert M_{k,n}\vvvert \vvvert M_{n,N}\vvvert.\)
Finally, to shorten the notations, we often omit the dependence in $\omega$, writing for example \(m_{k,n}(x)=m_{k,n}(x,\omega),\) and \(\vvvert M_{k,n}\vvvert=\Vert m_{k,n} \Vert_{\infty}= \sup_{x\in\XX} m_{k,n}(x,\omega).\)

\subsection{Assumptions}
We list here several hypotheses that will be used in the rest of the article.
\begin{ass} \label{ass: ergodicity}
The dynamical system $(\Omega, \mathcal{A}, \PP,\theta)$ is ergodic.
\end{ass}
We recall that a dynamical system is ergodic when any measurable set $A\in \mathcal{A}$ such that $\theta^{-1}(A)=A$ satisfies $\PP(A)\in \{0,1\}.$
\begin{ass} \label{ass: boundedness}
For almost all $\omega\in\Omega$, the function $x\mapsto m_{0,1}(x,\omega)$ is a positive function.
\end{ass}
By stationarity, \autoref{ass: boundedness} implies that $\PP(d\omega)$-almost surely, the product $M_{k,n}(\omega)$ is a continuous, non zero, positive linear operator.
We introduce the integrability property
\begin{ass} \label{ass: moments_m}
$\EE \left[ \log^+\Vert m_{0,1} \Vert_\infty \right]<\infty.$
\end{ass}
In particular, recalling that $\Vert m_{k,n} \Vert_\infty= \vvvert M_{k,n} \vvvert$, by submultiplicativity of the norm $\vvvert \cdot \vvvert$, \autoref{ass: ergodicity} and \autoref{ass: moments_m} imply that $\EE \log^+ \left(\left\vvvert M_{k,n} \right\vvvert\right)<\infty$ for all $k\leq n$.

We call admissible coupling constants a measurable map $$(\nu,c,d):\omega\in\Omega\mapsto (\nu_\omega,c(\omega), d(\omega))\in\mathcal{M}_1(\XX)\times[0,1]^2 $$ such that for $\PP$-almost any $\omega\in\Omega$,
\begin{itemize}
    \item[i)]for all $x\in\XX$ and all $f\in \B_+(\XX)$, the couple $(\nu_\omega,c(\omega))$ satisfies \begin{equation}\label{eq:prop c} \delta_x M(\omega)(f)\geq c(\omega) \Vert \delta_x M(\omega)\Vert \nu_\omega(f)\end{equation}
    \item[ii)] for all $n\geq 0$, the couple $(\nu_\omega,d(\omega))$ satisfies 
        \begin{equation} \nu_\omega (m_{1,n}) \geq d(\omega) \vvvert M_{1,n} \vvvert \label{eq: prop d}\end{equation}
\end{itemize}
When an admissible triplet is defined, we define the random variable \[ \gamma(\omega)=c(\omega)d(\omega).\]
Note that taking $c=d=0$ and any measurable map $\omega\mapsto \nu_\omega$ defines an admissible triplet, however in this case $\gamma=0$. 
{\color{black} Our main assumption is therefore
\begin{ass}
\label{ass : moments_cd} There exists an admissible triplet $\omega\mapsto(\nu_\omega,c(\omega),d(\omega))$ such that $\PP[\gamma>0]>0$.
\end{ass}}

%% file: 2_results.tex
\section{Statement of the results and structure of the paper}
\subsection{Main results}
Set \[ \bar{\gamma}:=\exp\left(\EE\left[ \log(1-\gamma)\right] \right)\in[0,1],\] and notice that \autoref{ass : moments_cd} yields $\bar{\gamma}<1$.
Under the previous assumptions, we prove the following Theorem:
\begin{theorem} \label{thme: approx}
Let $M:\Omega\longrightarrow \MX$ be a measurable map and assume that Assumptions \autoref{ass: ergodicity},\autoref{ass: boundedness} and \autoref{ass : moments_cd} hold. Then, 
\begin{enumerate}
\item[i)]$\PP(d\omega)$-almost surely, there exists a random function $h\in \linf$ such that, for any $\delta \in ({\color{black}\bar{\gamma}^{\PP[\gamma>0]}},1)$, for $n$ large enough, for any finite measures $\mu_1,\mu_2 \in \M_+(\XX)-\{0\}$, 
\begin{equation}\label{eq:approx}
    \left \Vert \mu_1 M_{0,n} - \frac{\mu_1(h)}{\mu_2(h)} \mu_2 M_{0,n} \right\Vert_{TV}\leq \delta^{n} \Vert \mu_1 M_{0,n}\Vert.
\end{equation}
Such a function $h$ is unique up to a multiplicative constant.
\item[ii)] There exists a probability measure $\Lambda$ on the space $\mathcal{M}_1(\XX)$, such that for any probability measure $\mu$, the sequence of random probability measures $\left(\mu \cdot M_{0,n}\right)$ converges in distribution towards $\Lambda$, in the space $\mathcal{M}_1(\XX)$, endowed with the total variation norm.
\item[iii)] Assuming additionally \autoref{ass: moments_m}, for almost any $\omega\in \Omega$ and any finite, positive, non-zero measure $\mu$,
\begin{equation} \label{eq : growth rate}
    \frac{1}{n}\log \Vert \mu M_{0,n} \Vert \underset{n\rightarrow\infty}{\longrightarrow}\inf_{N\in\NN}\frac{1}{N}\EE \left[ \log \vvvert M_{0,N}\vvvert \right]=\lambda\in[-\infty,\infty).
\end{equation}
\end{enumerate}
\end{theorem}

Note that the estimate \eqref{eq: presentation resultat} can be derived from Theorem \ref{thme: approx} by a choice of an arbitrary measure $\mu_2$, and by setting \[\pi_n=\mu_2\cdot M_{0,n}, \,  \, r_n=\frac{\Vert \mu_2 M_{0,n} \Vert}{\mu_2(h)},\]

The rest of this paper focuses on the independent case, that is, the case where the sequence of operators $(M_n)$ is i.i.d, with a law called $\P$. This can be obtained by setting $(\Omega,\A, \PP)$ to be the product space $\Omega=\left(\MX\right)^\NN$, $\PP=\P^{\otimes \NN}$ and $M:{\left(\MX\right)}^\NN \rightarrow \MX, (N_k)_{k\geq 0} \mapsto N_0$. In this independent case, we are able to characterize the measure $\Lambda$ on the projective space $\M_1(\XX)$ as the only invariant measure left invariant by the projective action of the matrices $(M_n)$. In other words, it is the only invariant probability measure of the Markov chain $(\mu_n)_{n\geq 0}$ defined by 
\[\mu_n= \frac{\mu_0 M_{0,n}}{\Vert \mu_0 M_{0,n} \Vert}\in\mathcal{M}_1(\XX).\]
A classical cocyle property yields the decomposition
\begin{equation}\label{eq:def cocycle}\log \Vert \mu_0 M_{0,n} \Vert_{TV}=\sum_{k=0}^{n-1} \log(\Vert \mu_k M_k \Vert)=\sum_{k=0}^{n-1} \rho(\mu_k,M_k),\end{equation}
where $(\mu_k,M_k)$ is a Markov chain on $\M_1(\XX)\times \MX$, with the unique invariant probability measure $\Lambda \otimes \P$.
These property allows us to obtain additional insight over the asymptotic behavior of the mass of the measure $\mu M_{0,n}$. To do so, we introduce the following strengthened versions of assumptions \autoref{ass: moments_m} and \autoref{ass : moments_cd}:
\begin{manualass}{\begin{NoHyper}\ref{ass: moments_m}\end{NoHyper}+} \label{ass: moments_m reinforced}
$\EE \left[ \left \vert \log \Vert m_{0,1} \Vert_\infty \right \vert \right]<\infty.$
\end{manualass}

\begin{manualass}{\begin{NoHyper}\ref{ass : moments_cd}\end{NoHyper}+} \label{ass: moments_cd reinforced} There exists an admissible triplet $(\nu,c,d)$ such that {$\EE \left| \log(\gamma)\right| <\infty$.}
\end{manualass}
Note that \autoref{ass: moments_cd reinforced} implies in particular that $\PP[\gamma>0]=1$.
Under these assumptions we are able to link the Lyapunov exponent $\lambda$, which governs the exponential growth of the mass of the measure $\Vert \mu M_{k,n} \Vert$, with the asymptotic projective distribution $\Lambda$ on $\mathcal{M}_1(\XX)$.
\begin{theorem}\label{thme:lyap}
Consider an  i.i.d sequence $(M_n)$ of elements of $\MX$ with law $\P$, suppose assumptions \autoref{ass: ergodicity}, \autoref{ass: boundedness}, \autoref{ass: moments_m reinforced} and \autoref{ass: moments_cd reinforced} hold. Then, the almost sure convergence \eqref{eq : growth rate} also holds in $\mathcal{L}^1(\Lambda\otimes \PP)$, that is
\[ \int \left| \frac{1}{n} \log \Vert \mu M_{0,n}\Vert - \lambda \right| d\Lambda(\mu) d\P^{\otimes n} (M_0,\cdots, M_{n-1})\underset{n\rightarrow \infty}{\longrightarrow} 0.\]
As a consequence, \begin{equation}\label{eq: lypaunov integral} \lambda = \int \log \Vert \mu M \Vert d\Lambda (\mu) d\P(M).\end{equation}
\end{theorem}
Moreover, in the critical case where $\lambda=0$, the law of large numbers \eqref{eq : growth rate} is not enough to know whether $\liminf \Vert \mu M_{0,n} \Vert$ and $\limsup \Vert \mu M_{0,n} \Vert$ are $0$, a positive real number, or $+\infty$. Answering these questions is the objective of our last theorem. When $\lambda=0$, the sequence $(M_n)$ is i.i.d and $\mu_0\sim \Lambda$, the increments $\rho(\mu_k,M_k)$ of the sum \eqref{eq:def cocycle} are centered and form an ergodic sequence. By analogy with a centered random walk with i.i.d increments, we expect their sum $\log \Vert \mu_0 M_{0,n} \Vert$ to oscillate between $-\infty$ and $+\infty$. More formally, we use the theory of Markov random walks, and in particular \cite{alsmeyer_recurrence_2001}, which establishes this oscillation property exists, except in some case called Null-Homology. In our context, we say that there is Null Homology when there exists some function
 $\eta : \mathcal{M}_1(\XX)\mapsto \RR$ such that $d(\Lambda \otimes \P)(\mu,M)$-a.s.
\begin{equation} \tag{NH} \label{eq :NH}\log \Vert \mu M \Vert = \eta( \mu \cdot M ) - \eta(\mu).\end{equation}
\stepcounter{equation}
When \eqref{eq :NH} does not hold, using results from \cite{alsmeyer_recurrence_2001}, we indeed obtain the oscillation of $\log \Vert \mu M_{0,n} \Vert$ for any initial measure $\mu$. On the contrary, when \eqref{eq :NH} holds, it is clear that $\log \Vert \mu_0 M_{0,n} \Vert= \eta(\mu_n)-\eta(\mu_0)$, thus $\log \Vert \mu_0 M_{0,n} \Vert$ may or may not oscillate, depending on $\eta$. More formally, our result is 
\stepcounter{equation}
\begin{theorem} \label{thme:osc}
Consider an  i.i.d sequence $(M_n)$ of elements of $\MX$ with law $\P$, suppose assumptions \autoref{ass: ergodicity}, \autoref{ass: boundedness}, \autoref{ass: moments_m} and \autoref{ass: moments_cd reinforced} hold. Assume additionally that $\lambda=0$. Then, if \eqref{eq :NH} does not hold, we have
\begin{equation} \tag{OSC}
\underset{n\rightarrow \infty}\liminf \log \Vert \mu M_{0,n} \Vert =-\infty \text{ and } \underset{n\rightarrow \infty}\limsup  \log \Vert \mu M_{0,n} \Vert =+\infty
\end{equation}
\addtocounter{equation}{1}
for any $\mu \in \mathcal{M}_+(\XX)-\{0\}.$
If \eqref{eq :NH} holds, let $\mu_0\sim \Lambda$, and let us note $a<b\in [-\infty,+\infty]$ the respective infinimum and supremum of the support of the random variable $\eta(\mu_0)$. It holds 
\[a=-\infty \Leftrightarrow \PP-\text{a.s., for any }\mu\in \mathcal{M}_+(\XX)-\{0\}, \, \underset{n\rightarrow \infty}\liminf \log \Vert \mu M_{0,n}\Vert=-\infty\]
and 
\[b=+\infty \Leftrightarrow \PP-\text{a.s., for any }\mu\in \mathcal{M}_+(\XX)-\{0\}, \, \underset{n\rightarrow \infty}\limsup \log \Vert \mu M_{0,n}\Vert=+\infty.\]
\end{theorem}
The notion of Null Homology already appears in the framework of products of $p\times p$ non-negative matrices, and \cite{hennion_limit_1997} provides some geometric condition which prevents Null-Homology (see in particular Theorem 5 of \cite{hennion_limit_1997}). We are unfortunately not able to generalize this condition in the infinite dimensional case.

\subsection{Structure of the paper}
Section \ref{section:proofs} contains the proofs of Theorems \ref{thme: approx}, \ref{thme:lyap} and \ref{thme:osc}. More precisely, in Subsection \ref{subs:contraction}, we recall how the coefficient $\gamma$ allow to control some contraction rates of the operators $M_{0,n}$. These results are adapted from \cite{bansaye_ergodic_2020}. In Subsection \ref{subs:asymptotic estimates}, we use the ergodic structure, in particular Assumptions \autoref{ass: ergodicity} and \autoref{ass : moments_cd} to obtain a geometric decay of the error terms that appeared in our previous estimations. In Subsection \ref{subs:thme} we derive the three claims of Theorem \ref{thme: approx}. Finally, in Subsection \ref{sec: independent}, we focus on the case where the sequence $(M_n)$ is i.i.d. In this case, a study of the invariant measures and the ergodicity properties of the Markov chains $(\mu_0 \cdot M_{0,n})_{n\geq 0}$ and $(\mu_0 \cdot M_{0,n}, M_n)_{n\geq 0}$, allows to prove Theorems \ref{thme:lyap} and \ref{thme:osc}.

Section \ref{sec: Hilbert} is dedicated to a comparison of our results with those obtained based on Hilbert contractions. More precisely, we show how natural conditions coming from Hilbert contractions techniques provide more tractable sufficient conditions for our Assumptions (in particular \autoref{ass : moments_cd}), both in finite and infinite dimension.

In Section \ref{sec: branching models}, we apply our results to study products of infinite Leslie matrices. This constitutes an example of an interesting class of systems that cannot be studied using the Hilbert metric. More precisely, we provide in Subsection \ref{subs: assumptions_age} reasonable sufficient conditions under which a product of Leslie matrices modelling the behavior of an age structured population satisfies assumptions \autoref{ass: ergodicity} to \autoref{ass : moments_cd}. However, when these conditions are not satisfied, it can be quite difficult to exhibit an admissible triplet such that $\gamma>0$, even on a deterministic and constant sequence of Leslie matrices.  To illustrate this fact, we present in Subsection \ref{subs : contrex} an example of system where $\gamma=0$ even if all the other assumptions are satisfied.

%% file: 3_proofs.tex
\section{Proofs}\label{section:proofs}
\subsection{Contraction results based on an inhomogeneous Doeblin minoration} \label{subs:contraction}
Given an admissible triplet $(\nu,c,d)$, we define for every $k\geq 0$ the $[0,1]$-valued random variables $c_k=c_k(\omega)=c(\theta^k(\omega))$, $d_{k+1}=d_{k+1}(\omega)=d(\theta^k(\omega))$ and $\gamma_k=\gamma_k(\omega)=\gamma(\theta^k(\omega))$, as well as the $\M_1(\XX)$ valued random variable $\nu_k=\nu_{\theta^k(\omega)}.$ 
Notice that $c=c_0$ and $d=d_1$, that the sequences $(c_k)_{k\geq 0}$, $(d_k)_{k\geq 1}$, $(\gamma_k)_{k\geq 0}$ are all stationary sequences and that for all $k\geq 0$, $$\gamma_k=c_kd_{k+1}.$$
Moreover the random variables $\nu_k,c_k,d_{k+1}$ satisfy some time-shifted versions of \eqref{eq:prop c} and \eqref{eq: prop d}. Indeed, for any $f\in\mathcal{B}_+(\XX)$ and any $x\in\XX$, it holds
\begin{equation}\label{eq:prop c shifted } \delta_x M_{k,k+1}(f)\geq c_k m_{k,k+1}(x) \nu_k(f),\, \PP\text{-a.s.}\end{equation}
and for any $n\geq k$,
\begin{equation} \nu_k (m_{k+1,n}) \geq d_{k+1} \vvvert M_{k+1,n} \vvvert = d_{k+1} \sup_{x\in\XX} m_{k+1,n}(x),\, \PP\text{-a.s.}\label{eq: prop d shifted}\end{equation}
The first step towards proving Theorem \ref{thme: approx} is establishing
\begin{proposition} \label{prop: contraction}
Suppose \autoref{ass: boundedness} and let $\omega\mapsto (\nu_\omega,c(\omega),d(\omega))$ be an admissible triplet. Then, $\PP(d\omega)$-almost surely, for any $k\leq n\leq N$ and any finite measures $\mu_1,\mu_2\in \mathcal{M}_+(\XX)-\{0\}$, it holds 
\begin{equation} \left\Vert \mu_1 \cdot M_{k,n}-\mu_2 \cdot M_{k,n}\right\Vert_{TV}\leq 2 \prod_{i=k}^{n-1}(1-\gamma_i), \label{eq: contraction distrib types} \end{equation}
and, if $n\geq 1$, it also holds $\PP(d\omega)$-a.s.
 \begin{equation} {\gamma_{n-1}}\left\vert \frac{\mu_1(m_{k,N})}{\mu_1(m_{k,n})}-\frac{\mu_2(m_{k,N})}{\mu_2(m_{k,n})}\right\vert\leq{2}\frac{\mu_2(m_{k,N})}{\mu_2(m_{k,n})} \prod_{i=k}^{n-1}(1-\gamma_i) .  \label{eq:contraction taux de croissance}\end{equation}
\end{proposition}

This result was already introduced in \cite{bansaye_ergodic_2020} in a somewhat different setup. We have chosen to state and prove it here for the sake of completeness. Its proof is based on performing a Doeblin minoration on a well-chosen sequence of auxiliary Markov operators $(P_{k,n}^N)$. This Doeblin property yields \eqref{eq: contraction distrib types}, a contraction property for the projective action of $M_{k,n}$ on the space of measures $\mathcal{M}_+(\XX)$. We derive then Equation \eqref{eq:contraction taux de croissance}, which describes how the growth of the mass $\Vert \mu M_{k,t}\Vert $ between times $t=n$ and $t=N$ depends on the initial measure $\mu$.

Let us introduce now the operators $P_{k,n}^N$ upon which we perform the desired Doeblin minoration. Under assumption \autoref{ass: boundedness}, $\PP(d\omega)$-almost surely, for any $k\leq n \leq N$, the functions $x\mapsto m_{k,n}(x)$ and  $x\mapsto m_{n,N}(x)$ are positive on $\XX$. For each $k\leq n\leq N$, an operator $P_{k,n}^N(\omega)$ can be defined $\PP(d\omega)$-almost surely, as follows : for each $x\in \XX$, for each positive measurable $f:\XX\longrightarrow\RR $, 
\[\delta_x P_{k,n}^N f = \frac{ \delta_x M_{k,n} (f m_{n,N})}{m_{k,N}(x)},\]
$P_{k,n}^N$ is a positive and conservative operator  (i.e. $P_{k,n}^N \mathds{1}=\mathds{1}).$ Indeed, by Equation \eqref{eq: semigroup}, for any $x\in\XX$, \[\delta_x P_{k,n}^N\mathds{1}=\frac{\delta_x M_{k,n} m_{n,N}}{\delta_x M_{k,N} \mathds{1}}=1.\]
Moreover, $P_{k,n}^N$ satisfies the relation :
\[P_{k,n}^N= P_{k,k+1}^N\cdots P_{n-1,n}^N.\]
Note that $P_{k,n}^N$ is a matrix when $\XX$ is countable and then, for any $x,y\in\XX$, \[ \mat{P_{k,n}^N}{x,y}=\frac{m_{n,N}(y)}{m_{k,N}(x)} \mat{M_{k,n}}{x,y}.\]

These operators satisfy a Doeblin contraction property summed up in 
\begin{lemma}\label{lem:Doeblin}
Assume Assumption \autoref{ass: boundedness} holds  and let $\omega\mapsto (\nu_\omega,c(\omega),d(\omega))$ be an admissible triplet. Then $\PP(d\omega)$-almost surely, all the $P_{k,n}^N(\omega)$ are well defined, and it holds
\begin{enumerate}
\item[i)] For any $n\leq N-1$,  there exists a random probability measure $\nu_{n,N}$ on $\XX$ such that, for any $x\in\XX$, \[\delta_x P_{n,n+1}^N\geq c_n d_{n+1} \nu_{n,N}=\gamma_n \nu_{n,N}.\]
\item[ii)] For any signed measures $\rho_1,\rho_2$, of same mass and any $n\leq N-1$,
\[\left \Vert \rho_1P_{n,n+1}^N -\rho_2 P_{n,n+1}^N \right\Vert_{TV}\leq (1-\gamma_n) \Vert\rho_1-\rho_2 \Vert_{TV}.\]
\item[iii)] For any $k\leq n\leq N$ and any signed measures $\rho_1,\rho_2$ of same mass,
\[\left \Vert \rho_1P_{k,n}^N-\rho_2P_{k,n}^N \right\Vert_{TV}\leq \prod_{i=k}^{n-1}(1-\gamma_i)\left \Vert \rho_1-\rho_2 \right\Vert_{TV}.\]
\end{enumerate}
\end{lemma}
Notice that in this lemma, our single assumption is \autoref{ass: boundedness}. It allows the $(P_{k,n}^N)$ to be defined $\PP(d\omega)$-almost surely. In particular, \autoref{ass : moments_cd} is not assumed, we allow $\gamma_n(\omega)=0$, in which case we just obtain that $P_{n,n+1}^N$ is $1$-contracting. 
\begin{proof}[Proof of Lemma \ref{lem:Doeblin}]
Let $\omega \in \Omega$ such that all the $P_{k,n}^N$ are well defined. For any $x\in\XX$ and any $f\in\mathcal{B}_+(\XX),$ it holds,  \[ \delta_x M_{n,n+1}(fm_{n+1,N})\geq c_n m_{n,n+1}(x) \nu_n(fm_{n+1,N}),\] thus \[\delta_x P_{n,n+1}^{N}f=\frac{\delta_x M_{n,n+1}(f m_{n+1,N})}{m_{n,N}(x)}\geq c_n\frac{\nu_n(fm_{n+1,N}) m_{n,n+1}(x)}{m_{n,N}(x)},\]
with, by definition of $d_{n+1}$ :
\[d_{n+1} m_{n,N}(x)=d_{n+1}\delta_x M_{n,n+1}(m_{n+1,N}) \leq \nu_n(m_{n+1,N})m_{n,n+1}{(x)}.\]
Therefore,
\[\delta_x P_{n,n+1}^{N}f\geq c_n d_{n+1}\frac{\nu(fm_{n+1,N})}{\nu(m_{n+1,N})}= c_n d_{n+1}\nu_{n,N}(f)=\gamma_n \nu_{n,N}(f),\]
    setting \[\nu_{n,N}(\cdot)=\frac{\nu_n(\cdot \, m_{n+1,N})}{\nu_n(m_{n+1,N})},\] which is a probability measure.
This concludes the proof of \textit{i)}.\\ Let us prove now \textit{ii)}. This result is a classical consequence of the previous point using the theory of Markov operators. A Markov operator $P$ is said to be $\delta$-Doeblin (with $\delta>0$) when there exists a probability measure $\mu$ such that $\delta_x P f \geq \delta \mu(f)$ for any $x$ in the state space and any $f\in\mathcal{B}_+(\XX)$. Furthermore, such an operator is $1-\delta$ contracting in total variation : for any signed measure $\rho$ of mass $0$, \[ \Vert \rho P \Vert_{TV}  \leq (1-\delta) \Vert \rho \Vert_{TV}.\] 
This property trivially holds for $\delta=0$ : any positive operator satisfies $\delta_x P f \geq 0$ when $f$ is a non negative function, and any Markov operator is $1=(1-0)$-contracting in total variation. 
In our context, the previous point of the lemma yields that $\PP(d\omega)$-almost surely, for any $n\leq N-1,$ and any $\nu\in \mathcal{M}_1(\XX)$, the Markov operator $P_{n,n+1}^N$ is $\gamma_n$-Doeblin. Therefore, for any $\rho_1,\rho_2 \in \lun$, such that $\rho_1(\mathds{1})=\rho_2(\mathds{1})$, noting $\rho=\rho_1-\rho_2$, it holds
\[ \left \Vert \rho P_{n,n+1}^N \right\Vert_{TV}\leq (1-\gamma_n) \Vert \rho\Vert_{TV}.\]
This proves \textit{ii)}, let us move now to \textit{iii)}. Since all the $P_{n,n+1}^{N}$ are conservative operators, the image of a measure of mass $0$ by $P_{n,n+1}^{N}$ is a measure of mass $0$. The equality \(P_{k,n}^N=P_{k,k+1}^N \cdots P_{n-1,n}^N, \)
yields  
\(\left \Vert \rho P_{k,n}^N \right\Vert_{TV}\leq (1-\gamma_{n-1})\Vert \rho P_{k,n-1}^N \Vert_{TV}.\) By induction, we deduce \[ \left \Vert \rho P_{k,n}^N \right\Vert_{TV}\leq \prod_{i=k}^{n-1}(1-\gamma_i)\left \Vert \rho \right\Vert_{TV}.\]
This concludes the proof. 
\end{proof}

\begin{proof}[Proof of Proposition \ref{prop: contraction}]
Let us prove first Inequality \eqref{eq: contraction distrib types}. Applying Lemma \ref{lem:Doeblin}, \textit{iii)} with $n=N$ and $\rho= \delta_x-\delta_y$, we get, $\PP(d\omega)$-almost surely
\[ \left\Vert \delta_x P_{k,n}^n -\delta_y P_{k,n}^n\right\Vert_{TV} \leq \prod_{i=k}^{n-1}(1-\gamma_i) \left\Vert \delta_x - \delta_y \right\Vert_{TV}\leq 2\prod_{i=k}^{n-1}(1-\gamma_i).\]
Hence, for any $f\in \linf$ and  $x,y\in \XX$,
\[  \left|  \frac{\delta_x M_{k,n}(f)}{m_{k,n}(x)} - \frac{\delta_y M_{k,n}(f) }{m_{k,n}(y)} \right| \leq 2 \left\Vert f \right\Vert_\infty \prod_{i=k}^{n-1}(1-\gamma_i).\]
Let $\mu_1$ and $\mu_2$ be two positive measures. The inequality 
\[  \left|  \delta_x M_{k,n}(f)  -  m_{k,n}(x) \frac{\delta_y M_{k,n}(f) }{m_{k,n}(y)} \right| \leq 2 m_{k,n}(x) \left \Vert f \right\Vert_\infty \prod_{i=k}^{n-1}(1-\gamma_i),\]
yields, after integrating with respect to $\mu_1(dx)$ :
\[  \left|  \mu_1 M_{k,n}(f)  -  \mu_1(m_{k,n}) \frac{\delta_y M_{k,n}(f) }{m_{k,n}(y)} \right| \leq 2 \mu_1(m_{k,n}) \left \Vert f \right\Vert_\infty \prod_{i=k}^{n-1}(1-\gamma_i),\]
so that
\[  \left|  \mu_1 \cdot M_{k,n} (f) - \frac{\delta_y M_{k,n}(f) }{m_{k,n}(y)} \right| \leq 2\left \Vert f \right\Vert_\infty \prod_{i=k}^{n-1}(1-\gamma_i).\]
Integrating now with respect to $\mu_2(dy)$, we obtain
\[  \left| \mu_1 \cdot M_{k,n} (f) - \mu_2 \cdot M_{k,n} (f) \right| \leq 2\left \Vert f \right\Vert_\infty \prod_{i=k}^{n-1}(1-\gamma_i),\]
and finally 
\[  \left\Vert  \mu_1 \cdot M_{k,n} - \mu_2 \cdot M_{k,n} \right\Vert_{TV} \leq 2 \prod_{i=k}^{n-1}(1-\gamma_i).\]
Let us move now to the proof of Inequality \eqref{eq:contraction taux de croissance}. Applying Inequality \eqref{eq: contraction distrib types} to the function $x\mapsto m_{n,N}(x)$, one gets 
\[  \left\vert \frac{\mu_1(m_{k,N})}{\mu_1(m_{k,n})}-\frac{\mu_2(m_{k,N})}{\mu_2(m_{k,n})}\right\vert\leq 2\left \Vert m_{n,N}\right\Vert_\infty \prod_{i=k}^{n-1}(1-\gamma_i)=2\left \vvvert M_{n,N}\right\vvvert\prod_{i=k}^{n-1}(1-\gamma_i).\]
By \eqref{eq:prop c shifted } and \eqref{eq: prop d shifted} it holds, \[\label{eq : controle d} {d_n}\left\vvvert M_{n,N}\right\vvvert \leq  \nu_{n-1}(m_{n,N}) \] and \[\label{eq : controle c} \mu_2(m_{k,N})= \mu_2 M_{k,n-1} M_{n-1,n} m_{n,N} \geq c_{n-1} \mu_2(m_{k,n}) \nu_{n-1}(m_{n,N}).\]
Combining these identities we obtain, $\PP(d\omega)$-almost surely,
\[  {c_{n-1}d_n }\left\vert \frac{\mu_1(m_{k,N})}{\mu_1(m_{k,n})}-\frac{\mu_2(m_{k,N})}{\mu_2(m_{k,n})}\right\vert\leq {2} \frac{\mu_2(m_{k,N})}{\mu_2(m_{k,n})} \prod_{i=k}^{n-1}(1-\gamma_i).\]
Thus taking an infimum in $\nu$, this yields, $\PP(d\omega)$-almost surely,
\[  {\gamma_{n-1}}\left\vert \frac{\mu_1(m_{k,N})}{\mu_1(m_{k,n})}-\frac{\mu_2(m_{k,N})}{\mu_2(m_{k,n})}\right\vert\leq 2 \frac{\mu_2(m_{k,N})}{\mu_2(m_{k,n})} \prod_{i=k}^{n-1}(1-\gamma_i).\]
This ends the proof.
\end{proof}
\subsection{Asymptotic estimates under ergodicity assumptions} \label{subs:asymptotic estimates}
For any $n$ such that $\gamma_{n-1}\neq 0$ and any $k\leq n$, we define
\begin{equation}
\label{eq:RHS}\Gamma_{k,n}= \frac{1}{\gamma_{n-1}} \prod_{i=k}^{n-1}{(1-\gamma_i)},
\end{equation}
and we set $\Gamma_{k,n}=+\infty$ when $\gamma_{n-1}=0$.
With these notations, Equation \eqref{eq:contraction taux de croissance} can be rewritten
\begin{equation}
\left\vert \frac{\mu_1(m_{k,N})}{\mu_1(m_{k,n})}-\frac{\mu_2(m_{k,N})}{\mu_2(m_{k,n})}\right\vert\leq 2 \Gamma_{k,n} \frac{\mu_2(m_{k,N})}{\mu_2(m_{k,n})}.
\end{equation}
In this subsection, we use the ergodicity Assumption \autoref{ass: ergodicity}, as well as Assumption \autoref{ass : moments_cd}, which provide a control on the sequence of Doeblin coefficients $(\gamma_n)$. In the following lemma, we use these assumptions to establish a geometric decay of  both $\prod_{i=k}^{n-1}(1-\gamma_i)$ and $\Gamma_{k,n}$ as $n\rightarrow \infty$. More precisely, \eqref{eq: limsup prod} provides the geometric decay of $\prod_{k=1}^{n-1}(1-\gamma_i)$ using Birkhoff's ergodic theorem. To derive the decay of $\Gamma_{k,n}$ as $n\rightarrow \infty$ we need additionally to avoid that $\gamma_{n-1}$ is either $0$ or too close to $0$. Under Assumption \autoref{ass: moments_cd reinforced}, the coefficients $(\gamma_{n})_{n\geq 0}$ are almost surely all nonzero, and \eqref{eq : limsup 1/g} provides a sufficient control for the geometric decay of $\Gamma_{k,n}$. 
{\color{black} Under the weaker assumption \autoref{ass : moments_cd}, it is possible however that $\gamma_{n-1}=0$ for some values of $n$. We can however focus on some random times at which $\gamma_{n-1}$ is greater that some predetermined level $\varepsilon>0$. We define those random times by
\begin{equation} \label{eq: def arr}
\arr{n}{\varepsilon}=\max\{i\leq n |\gamma_{i-1}\geq \varepsilon \}\cup\{-\infty\}
\end{equation}
for all $n\geq 0$ and all $\varepsilon>0$. \eqref{eq:prod arr} establishes that the subsequence of $(\Gamma_{k,n})_{n\geq 0}$ associated with the random times $(\arr{n}{\varepsilon})_{n\geq0 }$ still enjoys some slower geometric decay.}  

%\color{black}
\begin{lemma} \label{prop:error terms}
Assume that assumptions \autoref{ass: ergodicity}, \autoref{ass: boundedness} and \autoref{ass : moments_cd} hold. Then for any $k\geq 0 $, 
\begin{equation} \left(\prod_{i=k}^{k+n-1}(1-\gamma_i) \right)^\frac{1}{n}\underset{n\rightarrow \infty}{\longrightarrow} \bar{\gamma}<1, \PP(d\omega)-\text{almost surely} \label{eq: limsup prod}
\end{equation}
Moreover, for all $\varepsilon>0$ such that $\PP[\gamma\geq \varepsilon]>0$, it holds
{\color{black}\begin{equation}\label{eq:prod arr} \limsup_{n\rightarrow +\infty} \left(\Gamma_{k,\arr{n}{\varepsilon}} \right)^\frac{1}{n}\leq\bar{\gamma}^{\PP[\gamma\geq \varepsilon]}<1, \PP(d\omega)-\text{almost surely}.
\end{equation}}
Finally, assuming additionally \autoref{ass: moments_cd reinforced}, it holds 
\begin{equation} \label{eq : limsup 1/g}
\lim_{n\rightarrow \infty} \left(\frac{1}{\gamma_n}\right)^{\frac{1}{n}}=1, \PP(d\omega)-\text{almost surely}.
\end{equation}
\end{lemma}

\begin{proof}[Proof of Lemma \ref{prop:error terms}]
We recall first that by definition, for any $\omega \in \Omega$
\[(1-\gamma_i)(\omega)=(1-\gamma)\circ\theta^i(\omega).\]
Notice then that, by Assumption \autoref{ass : moments_cd}, for almost every $\omega\in \Omega$, $\gamma(\omega)\in (0,1]$ . Thus $\log(1-\gamma(\omega))\in [-\infty,0)$ and $\bar{\gamma}=\exp\left( \EE\left[ \log(1-\gamma(\omega)) \right] \right)\in [0,1).$

Thus for any $k$, 
\[ \log\left[\left( \prod_{i=k}^{k+n-1} (1-\gamma_i) \right) ^{\frac{1}{n}} \right]= \frac{1}{n} \sum_{i=0}^{n-1} \log (1-\gamma_k) \circ \theta^{i}.  \]
Since $\theta$ is an ergodic map, by Birkhoff's ergodic theorem, for any $k\geq 0$,
\[\frac{1}{n} \sum_{i=0}^{n-1} \log (1-\gamma_{k+i}) \underset{n\rightarrow \infty}{\longrightarrow}\EE\left[\log(1-\gamma)\right]= \log \bar{\gamma}, \,  \PP(d \omega)-\text{ almost surely}.\]
This yields \eqref{eq: limsup prod}.
\color{black}{Moreover, under \autoref{ass : moments_cd}, for $\varepsilon$ small enough, by Birkhoff's ergodic theorem, it holds almost surely
\[\lim_{n\rightarrow \infty} \frac{1}{n} \#\{ 1\leq k\leq n | \gamma_{k-1}\geq \varepsilon\}=\PP[\gamma\geq \varepsilon]>0.\]
However 
\[ \arr{n}{\varepsilon}\geq \#\{ 1\leq k\leq n | \gamma_{k-1}\geq \varepsilon\} \]
by definition of $\arr{n}{\varepsilon}.$
Thus, almost surely
\[\liminf_{n\rightarrow\infty} \frac{\arr{n}{\varepsilon}}{n}\geq \PP[\gamma\geq \varepsilon]>0,\]
in particular, $\lim_{n\rightarrow \infty} \arr{n}{\varepsilon}=+\infty$ almost surely.

As a consequence, for any $k\geq 0$, it holds $\PP$-a.s., for $n$ large enough
\begin{align*}
\left(\Gamma_{k,\arr{n}{\varepsilon}}\right)^\frac{1}{n} &=
\left(\frac{1}{\gamma_{\arr{n}{\varepsilon}-1}}\prod_{i=k}^{\arr{n}{\varepsilon}-1}(1-\gamma_i) \right)^\frac{1}{n}\leq \\ & \left(\frac{1}{\varepsilon}\right)^\frac{1}{n} \exp\left(
\frac{\arr{n}{\varepsilon}-k}{n} \frac{1}{\arr{n}{\varepsilon}-k}\sum_{i=k}^{\arr{n}{\varepsilon}-1} \log(1-\gamma_i)\right)
\end{align*}
where almost surely 
\[\lim_{n\rightarrow \infty}\frac{1}{\arr{n}{\varepsilon}-k}\sum_{i=k}^{\arr{n}{\varepsilon}-1} \log(1-\gamma_i)=\log(\bar{\gamma}).\]
Hence
\begin{align*}
\underset{n\rightarrow\infty}\limsup \left(\Gamma_{k,\arr{n}{\varepsilon}}\right)^\frac{1}{n}\leq \exp\left(\liminf_{n\rightarrow \infty}\left(
\frac{\arr{n}{\varepsilon}}{n}\right)\log(\bar{\gamma})\right) \leq \bar{\gamma}^{\PP\left[ \gamma\geq \varepsilon\right]}.
\end{align*}
This concludes the proof of \eqref{eq:prod arr}.} \color{black}{Let us move to the proof of Inequality \eqref{eq : limsup 1/g}. 
\\ Notice first that since $\gamma_n\leq 1$ for all $n$,\[\underset{n\rightarrow \infty}{\liminf} \left(\frac{1}{\gamma_n}\right)^\frac{1}{n}\geq 1.\]
Let us prove now the converse inequality.
Let us define for each $b >1$, \[Y_n(b)=\frac{-\log(\gamma_n)}{\log(b)}\geq 0,\]
and
\[N_b= \sum_{n\in\NN_0} \mathds{1}_{\left(1/{\gamma_n}\right)^{\frac{1}{n}} >b}=\sum_{n\in\NN_0} \mathds{1}_{Y_n(b)>n}.\]
For a given value of $b$, the sequence $(Y_n(b))_{n\in \NN_0}$ is stationary, thus \[\EE (N_b)= \sum_{n\in \NN_0} \PP\left[Y_n(b)>n \right]=\sum_{n\in \NN_0} \PP\left[Y_0(b)>n \right].\] It is a well known fact that for a nonnegative random variable $Y$, \[\EE[Y]<\infty \Leftrightarrow  \sum_{n\geq 0} \PP(Y>n) <\infty.\]
By Assumption \autoref{ass : moments_cd}, it holds $\EE\left[ Y_0(b) \right]<\infty$ and $\EE \left[ N_b \right] <\infty$ for all $b>1$. Therefore, $\PP(d\omega)$-almost surely, $N_b<\infty$, thus $\left(\frac{1}{\gamma_n}\right)^\frac{1}{n}\leq b$ for $n$ large enough and $\underset{n\rightarrow \infty}\limsup \left(\frac{1}{\gamma_n}\right)^\frac{1}{n}\leq 1$.
Finally, 
\[\underset{n\rightarrow \infty}{\lim} \left(\frac{1}{\gamma_n}\right)^{\frac{1}{n}} = 1, \, \PP(d\omega)-\text{almost surely.}\]}
\end{proof}

Putting the estimates from Lemma \ref{prop:error terms} together with Proposition \ref{prop: contraction}, we obtain 
\begin{proposition} \label{prop:approx}
Assume \autoref{ass: ergodicity}, \autoref{ass: boundedness}, \autoref{ass : moments_cd} hold. Let $\varepsilon>0$ such that $\PP[\gamma\geq \varepsilon]>0$. Then, $\PP(d\omega)$-almost surely, for any $k\in \NN_0$, there exists a bounded, non-negative measurable function $h_k$ such that, for $n$ large enough, for any $\mu_1,\mu_2 \in \M_+(\XX)-\{0\}$,
\begin{equation}\label{eq:PF1}
    \left \Vert \mu_1 M_{k,n} - \frac{\mu_1(h_k)}{\mu_2(h_k)} \mu_2 M_{k,n} \right\Vert_{TV}\leq {\color{black} {\Delta_{k,n}^\varepsilon}} \Vert \mu_1 M_{k,n}\Vert,
\end{equation}
{\color{black} where almost surely, for $n$ large enough
\[{\Delta_{k,n}^\varepsilon}:=\frac{8\Gamma_{k,\arr{n}{\varepsilon}}}{1-2\Gamma_{k,\arr{n}{\varepsilon}}}.\]
is well-defined and positive, and 
\begin{equation}\label{eq:geomlim prod arr}
\limsup_{n\rightarrow \infty} \left({\Delta_{k,n}^\varepsilon}\right)^\frac{1}{n}\leq\bar{\gamma}^{\PP[\gamma\geq \varepsilon]}<1.
\end{equation}}
Furthermore, $\PP(d\omega)$-almost surely, a function $h_k$ satisfying \eqref{eq:PF1} is unique up to a multiplicative constant.
\end{proposition}

\begin{proof}[Proof of Proposition \ref{prop:approx}]
{\color{black} We recall first that by \eqref{eq:prod arr}, it holds almost surely
\begin{equation}\lim_{n\rightarrow \infty} \Gamma_{k,\arr{n}{\varepsilon}}=0.\label{eq:lim prod arr}\end{equation}
As a consequence, almost surely, for $n$ large enough, $2\Gamma_{k,\arr{n}{\varepsilon}}<1$, thus ${\Delta_{k,n}^\varepsilon}$ is indeed well-defined and positive. 
We also derive \eqref{eq:geomlim prod arr} from \eqref{eq:prod arr} and \eqref{eq:lim prod arr}. As a consequence, $\lim_{n\rightarrow \infty} {\Delta_{k,n}^\varepsilon}=0$ almost surely.
}

Let us assume now that there exists a positive function $h_k$ satisfying Inequality \eqref{eq:PF1}. Then, if $x,y\in\XX$, setting $\mu_1=\delta_x,\mu_2=\delta_y$ and applying this inequality to the constant function $\mathds{1}$, we almost surely get, for $n$ large enough
\[\left|m_{k,n}(x)-\frac{h_k(x)}{h_k(y)}m_{k,n}(y)\right|\leq {\color{black}{\Delta_{k,n}^\varepsilon}} m_{k,n}(x)=\underset{n\rightarrow \infty}o(m_{k,n}(x)).\]
Thus \[\underset{n\rightarrow \infty}{\lim}\frac{h_k(x)}{h_k(y)}\frac{m_{k,n}(y)}{m_{k,n}(x)}=1,\]
which readily implies that \[\frac{h_k(x)}{h_k(y)}=\underset{n\rightarrow \infty}\lim \frac{m_{k,n}(x)}{m_{k,n}(y)}.\] This yields the unicity of $h_k$ up to a multiplicative constant, when it exists. Let us now prove the existence of a function $h_k$ satisfying \eqref{eq:PF1}.

By Inequality \eqref{eq:contraction taux de croissance}, with $\mu_1=\delta_x, \mu_2=\delta_y$, one gets, $\PP(d\omega)$-almost surely, for any $k\leq n \leq N$:
\begin{equation} \label{eq : taux croissance}\left\vert \frac{m_{k,N}(x)}{m_{k,N}(y)}- \frac{m_{k,n}(x)}{m_{k,n}(y)} \right \vert \leq 2\Gamma_{k,n}\frac{m_{k,n}(x)}{m_{k,n}(y)},  \end{equation}
where the right-hand-side term is infinite on the event $\{\gamma_{n-1}=0\}.$
Setting 
\[ \diam_{k,n}(x,y)= \sup_{N_1,N_2\geq n} \left\vert \frac{m_{k,N_1}(x)}{m_{k,N_1}(y)}- \frac{m_{k,N_2}(x)}{m_{k,N_2}(y)} \right \vert, \]
this yields, for any $x,y\in\XX$,
\[ \diam_{k,n}(x,y)\leq4\Gamma_{k,n}\frac{m_{k,n}(x)}{m_{k,n}(y)}.\]
Exchanging the roles of $x,y$, one gets :
\begin{equation} \min \left[\diam_{k,n}(x,y),\diam_{k,n}(y,x)\right] \leq 4\Gamma_{k,n}. \end{equation}
For all $k\geq 0$, $x,y\in\XX$, both the sequences $(\diam_{k,n}(x,y))_{n\geq k},(\diam_{k,n}(y,x))_{n\geq k}$ are non-increasing, and as a consequence from the definition of $\arr{n}{\varepsilon}$, $\arr{n}{\varepsilon}\leq n$, thus it holds 
{\color{black}\begin{align}
\min\left[ \diam_{k,n}(x,y),\diam_{k,n}(y,x)\right]&\leq \min\left[ \diam_{k,\arr{n}{\varepsilon}}(x,y),\diam_{k,\arr{n}{\varepsilon}}(y,x)\right]\nonumber\\
&\leq  \left(4\Gamma_{k,\arr{n}{\varepsilon}}\right)\underset{n\rightarrow \infty}\longrightarrow 0,   \label{eq:controle diam}
\end{align}}
$\PP(d\omega)$-almost surely, for $n$ large enough, by \eqref{eq:lim prod arr}.
Thus one of the sequences $$(\diam_{k,n}(x,y))_{n\geq k}, (\diam_{k,n}(y,x))_{n\geq k}$$ has $0$ as an adherence value. Since these sequences are non decreasing and positive, they converge, thus, $\PP$-a.s., one of them tends to $0$. 
Without loss of generality, suppose that 
\[ \diam_{k,n}(x,y)\underset{n\rightarrow \infty}{\longrightarrow} 0.\]
Then, the sequence of positive real numbers $\left(\frac{m_{k,n}(x)}{m_{k,n}(y)}\right)_{n\geq k}$ is a Cauchy sequence, it converges to a nonnegative limit $l_k(x,y).$ 
{\color{black} Moreover, as a consequence from the definition of $\arr{n}{\varepsilon}$, for all $k\leq n \leq N$, it holds $\arr{n}{\varepsilon}\leq n \leq N $ and thus $k \leq k \wedge \arr{n}{\varepsilon} \leq N$. Applying \eqref{eq : taux croissance} yields therefore
\begin{equation}\label{eq: taux croissance 2}
\left\vert \frac{m_{k,N}(x)}{m_{k,N}(y)}- \frac{m_{k,k\wedge\arr{n}{\varepsilon}}(x)}{m_{k,k\wedge\arr{n}{\varepsilon}}(y)} \right \vert \leq 2\Gamma_{k,k\wedge \arr{n}{\varepsilon}}\frac{m_{k,k\wedge \arr{n}{\varepsilon}}(x)}{m_{k,k\wedge \arr{n}{\varepsilon}}(y)}
.\end{equation}}
In particular, letting $N\rightarrow \infty$ in \eqref{eq: taux croissance 2} proves that the limit $l_k(x,y)=\lim_{N\rightarrow \infty}\frac{m_{k,N}(x)}{m_{k,N}(y)}$ satisfies
\begin{equation} \label{eq: equation positivity h_k} {\color{black}\left|l_k(x,y)-\frac{m_{k,k\wedge \arr{n}{\varepsilon}}(x)}{m_{k,k\wedge\arr{n}{\varepsilon}}(y)}\right| \leq 2\Gamma_{k,k\wedge \arr{n}{\varepsilon}}\frac{m_{k,k\wedge \arr{n}{\varepsilon}}(x)}{m_{k,k\wedge \arr{n}{\varepsilon}}(y)}}\end{equation}
Almost surely, it holds moreover for $n$ large enough
\[{\color{black} k\wedge \arr{n}{\varepsilon} =\arr{n}{\varepsilon} \text{ and }} 2\Gamma_{k,k\wedge \arr{n}{\varepsilon}}\leq \frac{1}{4}.\]
Plugging this into Equation \eqref{eq: equation positivity h_k} yields, for $n$ large enough
\[\left|\frac{m_{\color{black}k,\arr{n}{\varepsilon}}(x)}{m_{\color{black}k,\arr{n}{\varepsilon}}(y)}-l_k(x,y)\right|  \leq \frac{1}{4} \frac{m_{\color{black}k,\arr{n}{\varepsilon}}(x)}{m_{\color{black} k,\arr{n}{\varepsilon}}(y)}.\]
Since $\frac{m_{k,\arr{n}{\varepsilon}}(x)}{m_{k,\arr{n}{\varepsilon}}(y)}>0$, this implies that $l_k(x,y)=\lim_{n\rightarrow \infty } \frac{m_{k,n}(x)}{m_{k,n}(y)}>0$ and consequently, 
\[ \frac{m_{k,n}(y)}{m_{k,n}(x)} \underset{n\rightarrow\infty}{\longrightarrow} \frac{1}{l_k(x,y)}<\infty.\]
Note that Proposition \ref{prop: contraction} allows to prove that $\PP(d\omega)$-almost surely, \eqref{eq : taux croissance} holds jointly for any $k\leq n \leq N$ and any $x,y\in \XX$, thus so does \eqref{eq:controle diam}. Thus $\PP(d\omega)$-almost surely, all the sequences of the form $\left(\frac{m_{k,n}(x)}{m_{k,n}(y)}\right)_{n\geq k}$ for all $k\geq 0, x,y\in\XX$ converge to a positive limit as $n\longrightarrow \infty$.
\\ Now, let us fix an arbitrary element $x_0\in\XX$, and set $h_k(x)=\underset{n\rightarrow \infty}{\lim} \frac{m_{k,n}(x)}{m_{k,n}(x_0)}$ for all $x$. The function $h_k$ is positive and satisfies, for any $x,y\in\XX$,
$$\frac{m_{k,N}(x)}{m_{k,N}(y)}\underset{N\rightarrow \infty}\longrightarrow \frac{h_k(x)}{h_k(y)}.$$
Plugging this into Equation \eqref{eq: taux croissance 2}, we obtain almost surely, for $n$ large enough
{\color{black}
\begin{equation} \left \vert\frac{h_k(x)}{h_k(y)} - \frac{m_{k,\arr{n}{\varepsilon}}(x)}{m_{k,\arr{n}{\varepsilon}}(y)}\right \vert \leq 2\Gamma_{k,\arr{n}{\varepsilon}}\frac{m_{k,\arr{n}{\varepsilon}}(x)}{m_{k,\arr{n}{\varepsilon}}(y)}. \label{eq : cv h} \end{equation}}
Consequently,
\begin{align*} h_k(x) & \leq h_k(y) \left(1+2\Gamma_{k,\arr{n}{\varepsilon}}\right)\frac{m_{\color{black}k,\arr{n}{\varepsilon}}(x)}{m_{\color{black}k,\arr{n}{\varepsilon}}(y)}
\\ &\leq h_k(y)\left(1+2\Gamma_{k,\arr{n}{\varepsilon}}\right)\frac{\Vert m_{\color{black}k,\arr{n}{\varepsilon}}\Vert_{\infty}}{m_{\color{black}k,\arr{n}{\varepsilon}}(y)}<\infty,\end{align*}
by \autoref{ass: boundedness}, which implies that $h_k$ is bounded.
\\ {\color{black} Notice now that for $n$ large enough, taking $N=n$ in \eqref{eq: taux croissance 2} yields
\begin{equation}\label{eq: taux de croissance n et n tronque}
\left \vert \frac{m_{k,n}(x)}{m_{k,n}(y)}- \frac{m_{k,\arr{n}{\varepsilon}}(x)}{m_{k,\arr{n}{\varepsilon}}(y)} \right \vert \leq 2\Gamma_{k,\arr{n}{\varepsilon}}\frac{m_{k,\arr{n}{\varepsilon}}(x)}{m_{k,\arr{n}{\varepsilon}}(y)}
\end{equation}
Therefore, combining \eqref{eq: taux de croissance n et n tronque} with \eqref{eq : cv h}, we obtain
\begin{align}
\left\vert \frac{h_k(x)}{h_k(y)}- \frac{m_{k,n}(x)}{m_{k,n}(y)} \right \vert  & \leq \left \vert \frac{m_{k,n}(x)}{m_{k,n}(y)}- \frac{m_{k,\arr{n}{\varepsilon}}(x)}{m_{k,\arr{n}{\varepsilon}}(y)} \right \vert+\left \vert \frac{m_{k,N}(x)}{m_{k,N}(y)}- \frac{m_{k,\arr{n}{\varepsilon}}(x)}{m_{k,\arr{n}{\varepsilon}}(y)} \right \vert \nonumber \\
& \leq 4\Gamma_{k,\arr{n}{\varepsilon}}\frac{m_{k,\arr{n}{\varepsilon}}(x)}{m_{k,\arr{n}{\varepsilon}}(y)} \label{eq: cv h 2}
\end{align}
Once again, we recall that almost surely, for $n$ large enough
 \[ 2\Gamma_{k,\arr{n}{\varepsilon}}<1.\]
Therefore \eqref{eq: taux croissance 2} yields $\PP$-a.s., for $n$ large enough
\begin{equation}
\frac{m_{k,\arr{n}{\varepsilon}}(x)}{m_{k,\arr{n}{\varepsilon}}(y)}\leq \frac{1}{1-2\Gamma_{k,\arr{n}{\varepsilon}}}\frac{m_{k,n}(x)}{m_{k,n}(y)}.
\end{equation}
Plugging this into \eqref{eq: cv h 2} finally yields
\begin{equation}\left\vert \frac{h_k(x)}{h_k(y)}- \frac{m_{k,n}(x)}{m_{k,n}(y)} \right \vert  \leq \frac{4\Gamma_{k,\arr{n}{\varepsilon}}}{1-2\Gamma_{k,\arr{n}{\varepsilon}}}\frac{m_{k,n}(x)}{m_{k,n}(y)}= \frac{\Delta_{k,n}^\varepsilon}{2} \frac{m_{k,n}(x)}{m_{k,n}(y)} \label{eq : taux croissance 4}
\end{equation} }
\\Moreover, for any positive and finite measure $\mu_1\in \mathcal{M}_+(\XX)$, any $y\in\XX$, integrating \eqref{eq : taux croissance 4} with respect to $\mu_1(dx)$, one gets 
\[ \left|\frac{\mu_1(h_k)}{h_k(y)}-\frac{\mu_1(m_{k,n})}{m_{k,n}(y)}\right| \leq {\color{black}  \frac{\Delta_{k,n}^\varepsilon}{2}} \frac{ \mu_1(m_{k,n})}{m_{k,n}(y)},\]
Thus 
\[ \left|  m_{k,n}(y) \mu_1(h_k)-\mu_1(m_{k,n}) h_k(y)\right| \leq {\color{black}  \frac{\Delta_{k,n}^\varepsilon}{2}}  \mu_1(m_{k,n}) h_k(y). \]
Integrating with respect to any positive and finite measure $\mu_2(dy)$, this yields
\[\left|  \mu_2(m_{k,n}) \mu_1(h_k)-\mu_1(m_{k,n}) \mu_2(h_k)\right| \leq {\color{black}  \frac{\Delta_{k,n}^\varepsilon}{2}}  \mu_1(m_{k,n}) \mu_2(h_k), \]
and finally,
\begin{equation} \label{eq : relative error h} \left| \frac{\mu_1(h_k)}{\mu_2(h_k)}-\frac{\mu_1(m_{k,n})}{\mu_2(m_{k,n})}\right| \leq {\color{black}  \frac{\Delta_{k,n}^\varepsilon}{2} } \frac{ \mu_1(m_{k,n})}{\mu_2(m_{k,n})}.\end{equation}

Let us prove now that $h_k$ satisfies Inequality \eqref{eq:PF1}. It holds 
\begin{equation*}
  \begin{split}\left\Vert \mu_1 M_{k,n}- \frac{\mu_1(h_k)}{\mu_2(h_k)}  {\mu_2 M_{k,n}}\right\Vert_{TV}\leq & \left\Vert \mu_1 M_{k,n}-\mu_1(m_{k,n})\mu_2 \cdot M_{k,n}\right\Vert_{TV} \\ & + \left\Vert \mu_1(m_{k,n})\mu_2 \cdot M_{k,n} - \frac{\mu_1(h_k)}{\mu_2(h_k)} {\mu_2 M_{k,n}} \right\Vert_{TV}.\end{split}  
\end{equation*}
On the one hand, applying Inequality \eqref{eq: contraction distrib types}, one has, almost surely, for $n$ large enough
\begin{align*} 
\left\Vert \mu_1 M_{k,n}-\mu_1(m_{k,n})\mu_2 \cdot M_{k,n}\right\Vert_{TV} & \leq \mu_1(m_{k,n}) \left\Vert\mu_1 \cdot M_{k,n}-\mu_2 \cdot M_{k,n}\right\Vert_{TV} \\
& \leq 2\prod_{i=k}^{n-1}(1-\gamma_i)\mu_1(m_{k,n}) \\
 & \leq {\color{black} 2\prod_{i=k}^{\arr{n}{\varepsilon}-1}(1-\gamma_i)\mu_1(m_{k,n})} \\
& \leq {\color{black}{ \frac{\Delta_{k,n}^\varepsilon}{2} }\mu_1(m_{k,n})} .\end{align*}
On the other hand, by Equation \eqref{eq : relative error h}, it holds $\PP$-a.s., for $n$ large enough
\[\begin{split} \left\Vert \mu_1(m_{k,n})\mu_2 \cdot M_{k,n} - \frac{\mu_1(h_k)}{\mu_2(h_k)} {\mu_2 M_{k,n}} \right\Vert_{TV} &\leq \mu_2(m_{k,n}) \left| \frac{\mu_1(m_{k,n})}{\mu_2(m_{k,n})}-\frac{\mu_1(h_k)}{\mu_2(h_k)}\right| \\ & \leq {\color{black}{ \frac{\Delta_{k,n}^\varepsilon}{2} }} { \mu_1(m_{k,n})}.\end{split}\]
Finally, for $n$ large enough, it indeed holds almost surely for $n$ large enough
\[\begin{split}\left\Vert \mu_1 M_{k,n}- \frac{\mu_1(h_k)}{\mu_2(h_k)}  {\mu_2 M_{k,n}}\right\Vert_{TV} &\leq {\color{black}{\Delta_{k,n}^\varepsilon}}\Vert \mu_1 M_{k,n}\Vert_{TV}.\end{split}\]
This ends the proof.
\end{proof}
The control provided in \eqref{eq:PF1} only holds for large enough values of $n$ which are random and depend on $k$. Under the stronger assumption \autoref{ass: moments_cd reinforced}, it is possible to obtain a control which is uniform in $k\leq n$.
\begin{proposition}\label{prop: ergo unif in time}
Assume \autoref{ass: ergodicity}, \autoref{ass: boundedness}, \autoref{ass: moments_cd reinforced} hold. Then almost surely, for any $k\leq n$ and any measure $\mu_1,\mu_2\in\mathcal{M}_+(\XX)-\{0\}$, it holds 
\begin{equation}
 \left \Vert \mu_1 M_{k,n} - \frac{\mu_1(h_k)}{\mu_2(h_k)} \mu_2 M_{k,n} \right\Vert_{TV}\leq 4 \Gamma_{k,n} \Vert \mu_1 M_{k,n}\Vert=\underset{n\rightarrow \infty}o(\delta^n\Vert \mu_1 M_{k,n}\Vert),
\end{equation}
for any $\delta\in(\bar{\gamma},1)$.
\end{proposition}
\begin{proof}[Proof of Proposition \ref{prop: ergo unif in time}]
We start by letting $N$ go to infinity in \eqref{eq : taux croissance} and obtain, almost surely, for all $k,n \in \NN_0$ and all $x,y\in\XX$
\begin{equation}\label{eq: taux croissance 10}
\left\vert \frac{h_k(x)(x)}{h_k(y)}- \frac{m_{k,n}(x)}{m_{k,n}(y)} \right \vert2\Gamma_{k,n}\frac{m_{k,n}(x)}{m_{k,n}(y)},
\end{equation}
where by \autoref{ass: moments_cd reinforced} almost surely, $\gamma_{n-1}>0$ for all $n\geq 1$. Replacing \eqref{eq : taux croissance 4} by \eqref{eq: taux croissance 10} in the proof of Proposition \ref{prop:approx}, we obtain first
\begin{equation} \label{eq : relative error h 2} \left| \frac{\mu_1(m_{k,n})}{\mu_2(m_{k,n})}-\frac{\mu_1(h_k)}{\mu_2(h_k)}\right| \leq 2\Gamma_{k,n}\frac{ \mu_1(m_{k,n})}{\mu_2(m_{k,n})}\end{equation}
for any non-zero positives measures $\mu_1,\mu_2$.
From this we derive both
\begin{align*}
\left\Vert \mu_1 M_{k,n}-\mu_1(m_{k,n})\mu_2 \cdot M_{k,n}\right\Vert_{TV}&\leq \mu_1(m_{k,n}) \left\Vert\mu_2 \cdot M_{k,n}-\mu_2 \cdot M_{k,n}\right\Vert_{TV}\\ & \leq 2\prod_{i=k}^{n-1}(1-\gamma_i)\mu_1(m_{k,n})   \\ &\leq 2\Gamma_{k,n}\mu_1(m_{k,n})
\end{align*}
and
\[\begin{split} \left\Vert \mu_1(m_{k,n})\mu_2 \cdot M_{k,n} - \frac{\mu_1(h_k)}{\mu_2(h_k)} {\mu_2 M_{k,n}} \right\Vert_{TV} &\leq \mu_2(m_{k,n}) \left| \frac{\mu_1(m_{k,n})}{\mu_2(m_{k,n})}-\frac{\mu_1(h_k)}{\mu_2(h_k)}\right| \\ & \leq 2\Gamma_{k,n} { \mu_1(m_{k,n})},\end{split}\]
which we combine to obtain
\[\begin{split}\left\Vert \mu_1 M_{k,n}- \frac{\mu_1(h_k)}{\mu_2(h_k)}  {\mu_2 M_{k,n}}\right\Vert_{TV}\leq &\left\Vert \mu_1 M_{k,n}-\mu_1(m_{k,n})\mu_2 \cdot M_{k,n}\right\Vert_{TV} \\ & + \left\Vert \mu_1(m_{k,n})\mu_2 \cdot M_{k,n} - \frac{\mu_1(h_k)}{\mu_2(h_k)} {\mu_2 M_{k,n}} \right\Vert_{TV} \\ &\leq 4\Gamma_{k,n}\Vert \mu_1 M_{k,n}\Vert.\end{split}\]

\end{proof}
\subsection{Proof of Theorem  \ref{thme: approx}}\label{subs:thme}
\begin{proof}[Proof of assertion i), Uniform geometric ergodicity]
Let us take $k=0$ in Proposition \ref{prop:approx}. Then, for any $\varepsilon>0$, $\PP(d\omega)$-almost surely, noting $h=h_0$, it holds for any finite and positive measures $\mu_1,\mu_2$, on $\XX$,
\[\left\Vert \mu_1 M_{0,n}- \frac{\mu_1(h)}{\mu_2(h)}  {\mu_2 M_{0,n}}\right\Vert_{TV}\leq {\color{black}{\Delta_{0,n}^\varepsilon}} \Vert \mu_1 M_{0,n}\Vert,\]
{\color{black}where, $\PP(d\omega)$-almost surely,
\[\limsup_{n\rightarrow \infty} \left( {\Delta_{0,n}^\varepsilon}\right)^\frac{1}{n} \leq \bar{\gamma}^{\PP[\gamma\geq \varepsilon]}\in [0,1).\] 
 Let now $\delta\in (\bar{\gamma}^{\PP[\gamma>0]},1)$. For $\varepsilon>0$ small enough, $\delta>\bar{\gamma}^{\PP[\gamma\geq \varepsilon]}$. As a consequence, $\PP(d\omega)$-almost surely, for $n$ large enough, (depending on $\omega$),
\[ {\Delta_{0,n}^\varepsilon}\leq \delta^n.\]}
Thus, $\PP(d\omega)$-almost surely, for any $\delta\in (\bar{\gamma},1),$ for $n$ large enough and any positive and finite measures $\mu_1,\mu_2$,
\[ \left\Vert \mu_1 M_{0,n}- \frac{\mu_1(h)}{\mu_2(h)}  {\mu_2 M_{0,n}}\right\Vert_{TV}\leq\delta^n\Vert \mu_1 M_{0,n}\Vert_{TV}.\]
This proves Equation \eqref{eq:approx}. \end{proof}

\begin{proof}[Proof of assertion ii)]
The proof relies on a classical time-reversal technique, see e.g. \cite{cogburn_ergodic_1984, orey_markov_1991}, or \cite{hennion_limit_1997} for a version that is closer to our context. As stated in \cite[II.10.4, pp.239-241]{cornfeld_ergodic_1982}, the ergodic system $(\Omega, \mathcal{A},\PP,\theta)$ can be extended as an invertible ergodic system $(\overline{\Omega},\overline{\mathcal{A}}, \overline{\PP}, \overline{\theta})$, such that $\Omega\subset \overline{\Omega}$, $\overline{\theta}|_{\Omega}=\theta$, and $\overline{\theta}$ is a bijective, bimeasurable, measure preserving and ergodic mapping.
The definitions of $M_{k,n}$, $c_n,d_n,\nu_n,\gamma_n$ can be naturally extended to all $k\leq n$ in $\ZZ$, and one still has $c_n=c\circ \overline\theta^n$, $d_n=d\circ\overline\theta^n$, $\gamma_n=\gamma\circ \overline\theta^n$ for $n\in\ZZ$. Assumption \autoref{ass : moments_cd} implies that all the $(\gamma_n)_{n\in\ZZ}$ are almost surely positive and have $\log$-moments.

Therefore, Lemma \ref{lem:Doeblin} and Proposition \ref{prop: contraction} extend to indexes $k\leq n\leq N \in \ZZ$.

For nonnegative $ n\leq N$, for any positive measures $\mu_1,\mu_2$ on $\XX$, one has in particular 
\begin{equation}\label{eq:contraction reverse}\left \Vert \mu_1 \cdot M_{-n,0}-\mu_2 \cdot M_{-n,0} \right \Vert_{TV} \leq 2 \prod_{i=0}^{n-1}(1-\gamma_{-i-1}).\end{equation}
With $\mu_2=\mu_1 M_{-N,-n}$, this yields :
\[\left \Vert \mu_1 \cdot M_{-n,0}-\mu_1 \cdot M_{-N,0} \right \Vert_{TV}\leq 2\prod_{i=0}^{n-1}(1-\gamma_{-i-1}).\]

Noticing that $\overline{\theta}$ is now an ergodic automorphism of the measured space $\overline\Omega$, and applying Birkhoff-Khinchin Ergodic Theorem as stated in \cite[Theorem 1, p.11]{cornfeld_ergodic_1982}, one gets, for almost any $\omega\in \overline \Omega$
\[ \frac{1}{n} \sum_{i=0}^{n-1} \log (1-\gamma_{-i-1}) \circ \overline\theta ^{-i} \underset{n\rightarrow \infty}{\longrightarrow} \EE \left[\log(1-\gamma_{-1})\right]=\EE \left[\log(1-\gamma)\right].\]
Thus 
\[ \left( \prod_{i=0}^{n-1}(1-\gamma_{-i-1})\right)^\frac{1}{n} \underset{n\rightarrow \infty}{\longrightarrow} \exp\left(\EE \left[\log(1-\gamma)\right]\right)=\bar{\gamma} <1.\]
Therefore, almost surely, the sequence $\left(\mu_1 \cdot M_{-n,0}\right)_{n\in\NN_0}$ is a Cauchy sequence in the space $\mathcal{M}_1(\XX)$ of probabilities on $\XX$, endowed with the total variation norm. It thus converges almost surely to a random probability $\pi_{\mu_1}$ on $\XX$. For any finite, positive non-zero measures $\mu_1,\mu_2$, plugging  $\pi_{\mu_1},\pi_{\mu_2}$ into \eqref{eq:contraction reverse}, one proves that for almost any $\omega$, 
\[\pi_{\mu_1}=\pi_{\mu_2}.\]
Thus, there exists a random probability $\pi$, such that, almost surely, for any positive measure $\mu$
\[\lim_{n\rightarrow\infty}\left \Vert \mu \cdot M_{-n,0}-\pi\right\Vert_{TV}=0.\]
By stationarity of $\overline\theta$, \[\mu \cdot M_{-n,0}\overset{d}{=} \mu \cdot M_{0,n},\]
which proves that, noting $\Lambda$ the distribution of $\pi$, \[\mu \cdot M_{0,n} \underset{n\rightarrow\infty}{\overset{d}{\longrightarrow}} \Lambda.\]
\end{proof}
{\color{black} Before proving the last assertion of Theorem \ref{thme: approx}, we need to present the following lemma which establishes a sharp lower bound on some triangular inequality. 
\begin{lemma}\label{lem: controle gamma}
Assume \autoref{ass: boundedness} holds. Then, for any $0\leq k< n $ and any $\mu \in \mathcal{M}_+(\XX)$
\begin{equation}\label{eq: equivalence gamma} \gamma_k \Vert \mu M_{0,k+1}\Vert_{TV} \vvvert M_{k+1,n} \vvvert \leq \Vert \mu M_{0,n} \Vert_{TV} \leq \Vert \mu M_{0,k+1} \Vert_{TV} \vvvert M_{k+1,n} \vvvert.\end{equation}
\end{lemma}
\begin{proof}[Proof of Lemma \ref{lem: controle gamma}]
Let $\mu\in \mathcal{M}_+(\XX)-\{0\}$ and $0\leq k<n$. 
\[\mu M_{0,n}=(\mu M_{0,k}) M_{k} M_{k+1,n}\geq c_k \Vert  \mu M_{0,k+1} \Vert_{TV} \nu_k M_{k+1,n},\]
thus
\[ \Vert\mu M_{0,n} \Vert_{TV} = \mu M_{0,n} \mathds{1} \geq c_k \Vert \mu M_{0,k+1} \Vert (\nu_k M_{k+1,n})(\mathds{1})=c_k \Vert \mu M_{0,k+1} \Vert_{TV}  \Vert \nu_k M_{k+1,n} \Vert_{TV}.\]
By definition of $d_{k+1}$, it holds additionally
\[\Vert \nu_k M_{k+1,n} \Vert \geq d_{k+1} \vvvert M_{k+1,n}\vvvert.\]
Combining the two previous inequalities, we get
\[\Vert\mu M_{0,n} \Vert_{TV} \geq c_{k} d_{k+1} \Vert \mu M_{0,k+1} \Vert_{TV} \vvvert M_{k+1,n} \vvvert=\gamma_k\Vert \mu M_{0,k+1} \Vert_{TV} \vvvert M_{k+1,n} \vvvert\]
To obtain the second inequality of the lemma we simply write
\begin{align*}
\Vert \mu M_{0,n} \Vert &= \mu M_{0,n} \mathds{1}\\
&= \mu M_{0,k+1}M_{k+1,n}\mathds{1} \\
& \leq \sup_{x\in\XX} \delta_x M_{k+1,n} \mathds{1} \times \mu M_{0,k+1}\mathds{1} \\
&\leq \vvvert M_{k+1,n} \vvvert \Vert \mu M_{0,k+1}\Vert_{TV}.    
\end{align*}
\end{proof}}
\begin{proof}[Proof of assertion iii)]
We notice first that for any fixed integer, the $\PP-$almost sure convergence \[ n^{-1}\log \vvvert M_{0,n}\vvvert\underset{n\rightarrow \infty}{\longrightarrow}\lambda:=\inf_{N\geq 1}\frac{1}{N}\EE \left[ \log \vvvert M_{0,N}\vvvert \right] \] is a classical consequence of Kingsman's subbaditive ergodic theorem and the subadditivity property
\[\log \left\vvvert M_{0,n+p} \right \vvvert\leq \log \vvvert M_{0,n}\vvvert + \log \left \vvvert M_{n,n+p} \right\vvvert. \]
Note that applying this theorem requires Assumption \autoref{ass: moments_m} to ensure the integrability of $\log ^+ \vvvert M_{0,n} \vvvert$.
{\color{black} It is also classical that the same convergence holds for shifted sequences : for any $k\geq 0$,
\begin{equation} \label{eq: conclusion kingman} n^{-1}\log \vvvert M_{k,n}\vvvert\underset{n\rightarrow \infty}{\longrightarrow}\lambda\end{equation}
almost surely.
\\ Let now $\mu$ be a positive, finite measure on $\XX$. Let $k\geq 0$.
From Lemma \ref{lem: controle gamma}, it holds for $n\geq k$
\[\gamma_k \Vert \mu M_{0,k}\Vert_{TV} \vvvert M_{k,n} \vvvert \leq \Vert \mu M_{0,n} \Vert_{TV} \leq \Vert \mu M_{0,k} \Vert_{TV} \vvvert M_{k,n} \vvvert.\]
Thus
\[ \left| \frac{1}{n} \log \Vert \mu M_{0,n} \Vert_{TV}- \frac{1}{n} \log \vvvert M_{k,n} \vvvert \right| \leq \frac{1}{n} \left( \left| \log \Vert \mu M_{0,k} \Vert \right| +\left| \log \gamma_k \right|\right)\underset{n\rightarrow+\infty}{\longrightarrow } 0\]
almost surely on the event $\{\gamma_k>0\}$.
Combining this last estimate with \eqref{eq: conclusion kingman} proves that the convergence \[ \frac{1}{n} \log \Vert \mu M_{0,n}\Vert \underset{n\rightarrow \infty}\longrightarrow \lambda\] holds almost surely on the event $\{\gamma_k>0\}$, thus also on the event $\bigcup_{k\geq 0} \{\gamma_k>0\}$. Under Assumptions \autoref{ass: ergodicity} and \autoref{ass : moments_cd}, $\PP[\bigcup_{k\geq 0} \{\gamma_k>0\}]=1$ by Birkhoff's ergodic theorem,  thus the convergence \eqref{eq : growth rate} holds $\PP(d\omega)$-almost surely.}
\end{proof}
\subsection{The independent case : proof of Theorems \ref{thme:lyap} and \ref{thme:osc}}
\label{sec: independent}
Let us introduce the Markov chain $(\mu_n)_{n\geq 0}$ with state space $\mathcal{M}_1(\XX)$, defined by $\mu_{n+1}=\mu_n\cdot M_{n}=\mu_0\cdot M_{0,n+1}$.
The process $\left(\mu_n,M_n\right)_{n\geq 0}$, is then clearly also a Markov chain with state space $\mathcal{M}_1(\XX)\times \MX$ and transition kernel :
\[ Qf(\mu,M)=\int f(\mu\cdot M , N) d\P (N).\]
We denote $\PP_{\bar{\chi}}$ the law of the Markov chain $\left((\mu_n,M_n)\right)_{n\geq 0}$ when $(\mu_0,M_0)$ is distributed according to a measure $\bar{\chi}$ on $\mathcal{M}_1(\XX)\times \MX$.
Theorems \ref{thme:lyap} and \ref{thme:osc} rely on the study of the invariant measures and the ergodicity properties of the Markov chains $(\mu_n)_{n\geq 0}$ and $(\mu_n, M_n)_{n\geq 0}$. In particular, we show that the limit distribution $\Lambda$ of the Markov chain $(\mu_n)$ is its only invariant distribution. This is stated in the following proposition.
\begin{proposition}\label{prop:invariant measure}
Suppose $(M_n)_{n\geq 0}$ is an i.i.d sequence of elements of $\MX$ distributed according to $\P$, and satisfying the assumptions of Theorem \ref{thme: approx}. Then 
\begin{enumerate}[i)]
    \item For any initial distribution $\chi$ on $\mathcal{M}_1(\XX)$, the Markov chain $(\mu_n)_{n\geq 0}$ converges weakly to $\Lambda$.
    \item $\Lambda$ is the only invariant measure of the Markov chain $(\mu_n)_{n\geq 0}$.
    \item $\Lambda\otimes \P$ is the only invariant measure of the Markov chain $\left(\mu_n,M_n\right)_{n\geq 0}$.
\end{enumerate}
As a consequence, the dynamical systems associated with the Markov chains $(\mu_n)_{n\geq 0}$ and $(\mu_n,M_n)_{n\geq 0}$ are ergodic.
\end{proposition}
To prove Proposition \ref{prop:invariant measure}, we need to define the convolution operation $\star$ between probability measures on $\M_1(\MX)$ as follows :
For any $\Q_1,\Q_2\in \M_1(\MX)$, $\Q_1\star \Q_2$ is the law of $N_1 N_2$, where $(N_1,N_2)\sim \Q_1\otimes \Q_2$. We note, for any $\Q\in \M_1(\MX)$, $\Q^{\star n}$ the $n$-th convolution power of $\Q$. As an example, if $N_0,\dots, N_{n-1}$ are i.i.d with law $\Q$, $\Q^{\star n}$ is simply the distribution of $N_{0,n}=N_{0}\cdots N_{n-1}$. Given a probability distribution $\chi$ on $\M_1(\XX)$ and $\Q$ on $\MX$, we also note $\chi \stardot \Q$ the law of $\mu\cdot N$, where $(\mu,N)\sim \chi\otimes \Q$. 
These operations, previously defined in \cite{bougerol_products_1985} in a finite dimensional context, satisfy some elementary properties, summed up in the following lemma.
\begin{lemma}\label{lem:laws}
Let $\Q_1,\Q_2,\Q_3$ be probability measures on $\MX$ and $\chi$ be a probability measure on $\M_1(\XX)$, it holds
\begin{enumerate}[i)]
\item $(Q_1\star Q_2)\star Q_3=Q_1\star (Q_2\star \Q_3)$,
\item $(\chi \stardot \Q_1)\stardot \Q_2=\chi\stardot (\Q_1\star \Q_2)$,
\item For each $\Q\in \M_1(\MX)$, $\chi\mapsto\chi\stardot\Q$ is continuous with respect to the topology of convergence in law on $\M_1(\M_1(\XX))$.
\end{enumerate}
\end{lemma}
\begin{proof}[Proof of Lemma \ref{lem:laws}]
Consider $(N_1,N_2,N_3)\sim \Q_1 \otimes \Q_2\otimes \Q_3$. It holds
\[N_1N_2N_3=(N_1N_2)N_3=N_1(N_2N_3),\]
with $(N_1N_2)N_3\sim (Q_1\star Q_2)\star Q_3$ and $N_1(N_2N_3)\sim Q_1\star (Q_2\star \Q_3)$. This yields \textit{i)}.
\\ Let us prove now point \textit{ii)}. Consider $(\mu,N_1,N_2)\sim \chi \otimes \Q_1 \otimes \Q_2$. It holds \[\mu\cdot(N_1N_2)=(\mu\cdot N_1)\cdot N_2,\]
with $\mu\cdot(N_1N_2)\sim \chi\stardot (Q_1\star Q_2)$ and $(\mu\cdot N_1)\cdot N_2\sim (\chi \stardot Q_1)\stardot Q_2$. This yields \textit{ii)}.
\\ Let us move to the proof of \textit{iii)}. Consider a sequence of probability measures  $(\chi_n)$ on $\M_1(\XX)$, converging in distribution to $\chi$. Let us show that $(\chi_n\stardot \Q)_{n\geq 0}$ converges in distribution towards $\chi\stardot \Q$. 
Let $f$ be a continuous, bounded function on $\M_1(\XX)$, it holds :
\[ \int f(\mu) d(\chi_n\stardot \Q)(\mu)= \int \int f(\mu\cdot N) d\chi_n(\mu)d\Q(N)=\int \chi_n(g_N)d\Q(N), \]
where, for each $N\in \MX$, the function $g_N:\mu \mapsto f(\mu\cdot N)$ is continuous and bounded. Thus
\[ \chi_n(g_N)=\int f(\mu\cdot N)d\chi_n(\mu) \rightarrow \int f(\mu\cdot N )d\chi(\mu)=\chi(g_N).\]
This yields, by dominated convergence, as $n\rightarrow\infty$,
\[ \int \chi_n(g_N)d\Q(N)=\int f(\mu) d(\chi_n\stardot \Q)(\mu)\longrightarrow \int \chi(g_N)d\Q(N)=\int\int f(\mu\cdot N)d\chi(\mu)d\Q(N),\]
which implies \textit{iii)}.

\end{proof}

\begin{proof}[Proof of Proposition \ref{prop:invariant measure}]
 Let $f$ be a continuous and bounded function on $\M_1(\XX)$, it holds
 \[ \chi\stardot \P^{\star n} (f)=\int f(\mu) d(\chi \stardot \P^{\star n})(\mu)=\int\int f(\mu \cdot M_{0,n})d\chi(\mu) d\P^{\star n} (M_{0,n}).\]
 However, for any $\mu\in \M_1(\XX)$, Theorem \ref{thme: approx}, $ii)$ states that $(\delta_\mu \stardot \P^{\star n})_{n\geq 0}$ converges weakly towards $\Lambda$. Thus, for any $\mu\in \M_1(\XX)$, as $n\rightarrow \infty$
 \[ \int f(\mu \cdot M_{0,n}) d\P^{\star n } (M_{0,n})=(\delta_\mu \stardot \P^{\star n})(f) \underset{n\rightarrow \infty}\longrightarrow \Lambda(f).\]
By dominated convergence, this yields
\[ \chi \stardot \P^{\star n}(f)=\int\int f(\mu \cdot M_{0,n})d\chi(\mu) d\P^{\star n} (M_{0,n})\underset{n\rightarrow \infty}\longrightarrow \Lambda(f),\]
which proves the weak convergence 
\[ \chi\stardot \P^{\star n} \underset{n\rightarrow\infty}{\longrightarrow}\Lambda \] in the metric space $\M_1(\M_1(\XX))$, for any probability distribution $\chi$.
This proves \textit{i)}.
\\ Since, by Lemma \ref{lem:laws}, \textit{iii)}, the map $\mu \mapsto \mu \stardot \P$ is continuous, this proves that $\Lambda$ is one of its fixed points, namely :
\[\Lambda = \Lambda \stardot \P.\]
On the other hand, if $\chi\stardot \P=\chi$, the sequence $(\chi\stardot \P^{\star n})_{n\geq 0}$ is constant and converges to $\chi$. By unicity of the limit, it holds 
\[ \chi=\Lambda.\]
This proves that $\Lambda$ is the only invariant measure of the Markov chain $(\mu_n)$, i.e. \textit{ii)}.
\\ Let $(\mu_0,M_0)\sim \Lambda \otimes \P$. Then $\mu_1= \mu_0 \cdot M_0\sim \Lambda \stardot \P$, $M_1\sim \P$ and $M_1$ is independent of $\mu_0$, $M_0$ and thus $\mu_1$. Therefore $(\mu_1,M_1)\sim \Lambda \otimes \P$, and $\Lambda \otimes \P$ is thus an invariant measure of the Markov chain $(\mu_n,M_n)_{n\geq 0}$.

Conversely, consider now a probability measure $\bar{\chi}$ on $\M_{1}(\XX)\times \MX$, suppose it is an invariant measure of the Markov chain $\left((\mu_n,M_n)\right)_{n\geq 0}$. 
The definition of the transition kernel $Q$ implies that $\mu_1=\mu_0\cdot M_0$, $M_1\sim \P$ and $M_1$ is independent of $(\mu_0,M_0)$, and therefore $M_1$ is independent of $\mu_1$.
However the second term $(\mu_1,M_1)$ of the Markov chain is distributed according to $\bar{\chi} Q=\bar{\chi}$ by invariance.
Thus $\bar{\chi}$ is of the form $\bar{\chi}=\chi\otimes \P$. 

Additionally, if $(\mu_0,M_0)\sim \bar{\chi} =\chi\otimes \P$, then $\mu_1=\mu_0\cdot M_0\sim \chi\stardot \P$. But by invariance of $\bar{\chi}$, $\mu_1 \sim \chi$, thus
\[ \chi \stardot \P=\chi.\]
By Proposition \ref{prop:invariant measure}, this implies that $\chi=\Lambda$.
Finally, this proves that $\Lambda\otimes \P$ is the only invariant measure of the Markov chain $\left(\mu_n,M_n\right)_{n\geq 0}$.
By Corollary 5.12 of \cite{hairer_ergodic_2018}, since both the processes $(\mu_n)_{n\geq 0}$ and $\left(\mu_n,M_n\right)_{n\geq 0}$ are Markov chains with a unique invariant measure, they are both ergodic.
\end{proof}

One additional lemma is required before proving Theorem \ref{thme:lyap}.

{\color{black} \begin{lemma} \label{lem:0-1}
Let $(M_n)_{n\geq 0}$ be an i.i.d sequence of elements of $\MX$ with law $\P$, satisfying the assumptions \autoref{ass: ergodicity}, \autoref{ass: boundedness} and \autoref{ass: moments_cd reinforced}.
The sets $$\{ \forall \mu \in \mathcal{M}_1(\XX), \underset{n\rightarrow\infty}\limsup \log \Vert \mu{M_{0,n}} \Vert=+\infty \},$$ $$\{ \exists \mu \in \mathcal{M}_1(\XX), \underset{n\rightarrow\infty}\limsup \log \Vert \mu{M_{0,n}} \Vert=+\infty \}$$ and the event $$\{ \underset{n\rightarrow\infty}\limsup \log \vvvert {M_{1,n}} \vvvert=+\infty \}$$ coincide up to $\P^{\otimes \NN}$-negligible events. A similar statement holds replacing $\limsup$ by $-\liminf$ in the three events. 
\end{lemma}
\begin{proof}[Proof of Lemma \ref{lem:0-1}]
From Lemma \ref{lem: controle gamma}, with $k=0$, we get 
\[ \gamma_0 \Vert \mu M_{0} \Vert_{TV} \vvvert M_{1,n} \vvvert \leq \Vert \mu M_{0,n}\Vert_{TV} \leq \Vert \mu M_0 \Vert_{TV} \vvvert M_{1,n}\vvvert ,\]
where $\gamma_0\Vert \mu M_{0} \Vert_{TV}>0$ almost surely by \autoref{ass: boundedness} and \autoref{ass: moments_cd reinforced}. The lemma is a straightforward consequence of this inequality. 
\end{proof}}

Let us prove now Theorem \ref{thme:lyap}.
\begin{proof}[Proof of Theorem \ref{thme:lyap}, \textit{i)}]  
Let us notice first that when $\mu$ is a probability measure, $\rho:(\mu, M)\mapsto \log \Vert \mu M \Vert $ satisfies the cocycle property
\begin{align}  \rho(\mu_0,M_{0,n})&=\log \Vert \mu_0 M_{0,n}\Vert \nonumber \\
&=\sum_{k=0}^{n-1} \log \left(\frac{\Vert \mu_0 M_{0,k+1} \Vert}{\Vert \mu_0 M_{0,k}\Vert}\right)\nonumber \\
&=\sum_{k=0}^{n-1} \log \Vert (\mu_0\cdot M_{0,k}) M_{k} \Vert \nonumber\\
&=\sum_{k=0}^{n-1} \rho( \mu_k, M_k ).\label{eq:cocycle}\end{align}
From Equation \eqref{eq:cocycle}, we derive
\[\frac{1}{n} \log \Vert \mu_0 M_{0,n}\Vert = \frac{1}{n} \sum_{k=0}^{n-1} \rho( \mu_k, M_k ).\]
By Birkhoff's Ergodic Theorem, since $(\mu_k,M_k)$ is an ergodic Markov chain, with stationary distribution $\Lambda \otimes \P$, this quantity converges $\PP_{\Lambda \otimes \P}$-almost surely and in $\mathcal{L}^1(\PP_{\Lambda \otimes \P})$ towards $\int \rho d(\Lambda\otimes \P)$ provided $\rho$ is an $\mathcal{L}^1$ function with respect to $\Lambda\otimes \P$. Let us check now this integrability property. 
Let $(\mu_0,{(M_n)}_{n\geq 0}) \sim \Lambda\otimes \P^{\NN}$. Then, applying \eqref{eq: equivalence gamma} with $n=2$ yields
\[ \gamma_0 \Vert \mu_0 M_{0}\Vert \vvvert M_{1} \vvvert \leq \Vert \mu_0 M_{0,2} \Vert \leq \Vert \mu_0 M_{0} \Vert \vvvert M_{1} \vvvert.\]
Noting $\mu_1=\mu_0 \cdot M_0$, we get, since $\Vert \mu_0 M_0\Vert \neq 0$ almost surely, 
\[ \gamma_0 \vvvert M_1 \vvvert \leq \Vert \mu_1 M_{1} \Vert \leq \vvvert M_1 \vvvert,\]
thus 
\[ \vert \rho(\mu_1,M_1) \vert\leq \vert \log \Vert M_1 \Vert \vert + \vert \log(\gamma_0) \vert. \]
Note that $(\mu_1,M_1)\sim \Lambda \otimes \P,$ since by definition $\mu_1=\mu_0\cdot M_0$ and $(\mu_0, (M_n)_{n\in \NN} ) \sim \Lambda \otimes \P^{\otimes \NN}.$
Thus, under \autoref{ass: moments_m reinforced} and \autoref{ass: moments_cd reinforced}, it holds
\[ \EE\left[\vert \rho(\mu_1,M_1)\vert \right] =\int \vert \rho \vert d(\Lambda \otimes \P) \leq \EE \vert \log \Vert M_1 \Vert \vert +\EE \vert \log(\gamma_0) \vert <\infty.\]
This proves that $\rho$ is integrable with respect to $\Lambda \otimes \P$, thus the convergence $$n^{-1}\log(\Vert \mu_0 M_{0,n}\Vert)\underset{n\rightarrow \infty}\longrightarrow \int \rho \,d\Lambda\otimes \P$$ holds in $\mathcal{L}^1(\Lambda\otimes \PP)$ and $\PP_{\Lambda \otimes \PP}$-almost surely.
Since by Theorem \ref{thme: approx}, almost surely, for all $\mu$, it holds $n^{-1}\log \Vert \mu M_{0,n} \Vert \underset{n\rightarrow \infty}\longrightarrow\lambda$, by unicity of the almost sure limit,
\[ \lambda= \int \rho\, d(\Lambda \otimes \P).\]
\end{proof} 

\begin{proof}[Proof of Theorem \ref{thme:osc}, \textit{ii)}]
Note now $X_n= \rho(\mu_{n-1},M_{n-1})$, for $n\geq 1$. Then it holds, for $n\geq 0$, for any probability measure $\mu_0$ 
\[ (\mu_{n+1}, X_{n+1})=(\mu_n\cdot M_n, \rho(\mu_n,M_n)).\]
Thus, $(\mu_n, X_n)_{n\geq 0}$ is a Markov chain on $\mathcal{M}_1(\XX)\times \RR$ such that \[ \PP\left[ (\mu_{n+1}, X_{n+1})\in A\times B |(\mu_n,X_n)\right]= \int \mathds{1}_A( \mu_n\cdot M) \mathds{1}_{B}\left(\rho(\mu_n,M)\right) d\P(M). \]
Thus $S_n=\log \Vert \mu M_{0,n} \Vert = X_1 + \cdots + X_{n}$ is a Markov random walk associated with $(\mu_n,X_n)$, in the sense of \cite{alsmeyer_recurrence_2001}. Suppose that $\lambda =0$. By Theorem \ref{thme: approx}, it holds  \[ n^{-1}S_n= n^{-1}\log \Vert \mu M_{0,n}\Vert \underset{n\rightarrow \infty}{\longrightarrow} 0,\] $d\PP_{\Lambda \otimes \P^{\otimes \NN}}\left(\mu,(M_n)_{n\geq 0}\right)$-almost surely, thus in probability with respect to $\PP_{\Lambda \otimes \P^{\otimes \NN} }$. Since moreover, $(\mu_n)$ is an ergodic Markov chain, the assumptions of \cite{alsmeyer_recurrence_2001} are satisfied. If there exists a function $\eta$ such that $\PP_{\Lambda \otimes \P}$-almost surely, for $n\geq 1$,
\begin{equation} \label{eq:cobord} X_n= \eta(\mu_n)-\eta (\mu_{n-1}),\end{equation}
then taking $n=1$ shows that we are in the case of Null Homology \eqref{eq :NH}. In this case, it holds moreover
\[ \log \Vert \mu_0 M_{0,n} \Vert =X_1+\cdots +X_{n}=\eta(\mu_{n})-\eta(\mu_0).\]
Thus, almost surely, noting $a,b$ the respective infimum and supremum of the support of $\eta(\mu)$, when $\mu \sim \Lambda$, since the sequence $(\mu_{n})$ is a stationary and ergodic sequence with law $\Lambda$, it holds
\[ \underset{n\rightarrow\infty}\liminf\log \Vert \mu_0 M_{0,n} \Vert = a-\eta(\mu_0),\text{ and } \underset{n\rightarrow\infty}\limsup \log \Vert \mu_0 M_{0,n} \Vert = b-\eta(\mu_0).\]
Thus the almost sure finiteness of these quantities are respectively equivalent to the finiteness of $a$ and $b$. 
If Equation \eqref{eq:cobord} does not hold, then we are in the setup of Theorem 2 or 3 of \cite{alsmeyer_recurrence_2001}. These two Theorems imply that the Markov Random Walk $(S_n)$ oscillates : $\limsup S_n=+\infty$ and $\liminf S_n=-\infty$ $\PP_{\Lambda \otimes \P}-$almost surely. However, by Lemma \ref{lem:0-1}, \textit{ii)}, this implies that $\P^{\otimes \NN}$-almost surely, for every $\mu \in \mathcal{M}_1(\XX),$ \[\underset{n\rightarrow \infty}{\limsup} \log \Vert \mu M_{0,n} \Vert=-\underset{n\rightarrow \infty}{\liminf} \log \Vert \mu M_{0,n} \Vert=+\infty.\]
This concludes the proof.
\end{proof}

%% file: 4_hennion.tex
\section{Sufficient conditions under uniform positivity assumptions} \label{sec: Hilbert}
In the finite dimensional case $\XX=\{1,\dots, p\}$, that is when studying products of $p\times p$ matrices, similar (and actually, more complete) results are obtained in \cite{hennion_limit_1997}. They rely on the very mild assumption 
\begin{ass} \label{ass: hennion}
\(\PP\left[\bigcup_{k\in\NN} \left\{M_{0,k}\in \accentset{\circ}{\mathcal{S}}\right\}\right]=1,\)
\end{ass}
where $\accentset{\circ}{\mathcal{S}}$ refers to the set of $p\times p$ matrices with positive entries. We expect that this approach, based on Hilbert contractions, might be extended in infinite dimensional contexts. This will require to introduce the notion of uniformly positive operators to strengthen the notion of positive matrices, and state an infinite dimensional generalization of \autoref{ass: hennion}, as we explain in Subsection \ref{subs: inFiDi Hilbert}.

This section aims at comparing our assumptions both with $\autoref{ass: hennion}$, and its natural generalization in infinite dimension.

We did not success in proving that $\autoref{ass: hennion}$ alone is enough for our assumptions to hold. However we provide mild additional assumptions that, together with \autoref{ass: hennion}, constitute sufficient conditions for our assumptions (\autoref{ass: boundedness}, \autoref{ass: moments_m}, \autoref{ass : moments_cd}, \autoref{ass: moments_cd reinforced}) to hold, and thus for Theorems \ref{thme: approx}, \ref{thme:lyap}, \ref{thme:osc} to apply.

\subsection{The finite dimensional case}\label{subs:}
Let us focus in this subsection on the case where $\XX$ is finite, let us note $p=\vert \XX \vert$. Consider a stationary and ergodic sequence $(M_n)_{n\in\NN_0}$ of $p\times p$ matrices with nonnegative entries. Checking whether Assumptions \autoref{ass: ergodicity}, \autoref{ass: boundedness}, \autoref{ass: moments_m} are satisfied is quite straightforward, since these three assumptions only involve the law of the first matrix of the sequence. Let us see now how the additional Assumption \autoref{ass: hennion} can help exhibit an admissible triplet in order to check that Assumptions \autoref{ass : moments_cd} holds.

\begin{lemma}
\label{lem: pos d hennion}
Consider a random and stationary sequence of $p\times p$ matrices $M_n= \left(\mat{M_n}{x,y}\right)_{x,y\in\XX}$ with nonnegative entries, satisfying \autoref{ass: boundedness} and \autoref{ass: hennion}.
Then for any measurable map $\omega\in\Omega\mapsto \nu(\omega)\in \mathcal{M}_1(\XX)$, there exists a random variable $d:\omega\mapsto d(\omega)$ such that $\PP[d>0]=1$ and \eqref{eq: prop d} holds.
\end{lemma} 

\begin{proof}
The following decomposition holds : for any $1\leq k \leq n$, $x\in\XX$, $\omega\in\Omega$
\[ m_{1,n}(x)=\delta_x M_{1,n}\mathds{1}= \delta_x M_{1,k} m_{k,n}=\sum_{z\in\XX} \mat{M_{1,k}}{x,z} m_{k,n}(z).\]
Thus
\[ \nu_0 ( m_{1,n}) =\nu_0  M_{1,n} \mathds{1}=\sum_{y,z\in \XX} \nu_0 ({y}) \mat{M_{1,k}}{y,z} m_{k,n}(z),\]
where we note $\nu_0=\nu_\omega$.
The fact that $\nu_0 $ is a probability measure yields,
\begin{align*} m_{1,n}(x) &\leq \sup_{y,z\in \XX}\left( \frac{\mat{M_{1,k}}{x,z}}{\mat{M_{1,k}}{y,z}}\right) \sum_{y\in\XX}\nu_0 (y)  \mat{M_{1,k}}{y,z} m_{k,n}(z)\\ &\leq \sup_{y,z\in \XX} \frac{\mat{M_{1,k}}{x,z}}{\mat{M_{1,k}}{y,z}} \nu_0 (m_{1,n}),\end{align*}
 for any $1\leq k \leq n $, with the conventions $\frac{\mat{M_{1,k}}{x,z}}{\mat{M_{1,k}}{y,z}}=0$ as soon as $\mat{M_{1,k}}{x,z}=0$ and $\frac{\mat{M_{1,k}}{x,z}}{\mat{M_{1,k}}{y,z}}=\infty$ if $\mat{M_{1,k}}{x,z}\neq 0$ and $\mat{M_{1,k}}{y,z}=0$.
Thus, for any $n\geq k$,
\begin{equation} \label{eq : majoration hennion 1} \frac{\vvvert M_{1,n} \vvvert }{\nu_0 (m_{1,n})}= \frac{ \Vert m_{1,n} \Vert_\infty}{\nu_0 (m_{1,n})} \leq \sup_{x,y,z\in \XX} \frac{\mat{M_{1,k}}{x,z}}{\mat{M_{1,k}}{y,z}}.\end{equation}
This yields \begin{equation} \label{eq : minoration hennion 1}\inf_{x,y,z\in\XX} \frac{\mat{M_{1,k}}{y,z}}{\mat{M_{1,k}}{x,z}}\leq \inf_{n\geq k} \frac{ \nu_0 (m_{1,n})}{\Vert m_{1,n}\Vert_\infty},\end{equation}
and therefore
\[ \sup_{k\geq 1} \inf_{x,y,z\in\XX} \frac{\mat{M_{1,k}}{y,z}}{\mat{M_{1,k}}{x,z}}\leq \liminf_n \frac{\nu_0 (m_{1,n})}{\vvvert M_{1,n}\vvvert}.\]
By Assumption \autoref{ass: hennion}, $\PP(d\omega)$- almost surely, there exists a random integer $k_\omega$ such that $M_{1,k_\omega}(\omega)\in \accentset{\circ}{\mathcal{S}}$. Since $\XX$ is finite, we get, for $\PP$-almost any $\omega$,
\[1<\inf_{x,y,z\in\XX} \frac{\mat{M_{1,k_\omega}}{y,z}}{\mat{M_{1,k_\omega}}{x,z}}<\sup_{k\geq 1}\inf_{x,y,z\in\XX} \frac{\mat{M_{1,k}}{y,z}}{\mat{M_{1,k}}{x,z}}\leq \liminf_n \frac{\nu_0 (m_{1,n})}{\vvvert M_{1,n}\vvvert}.\]
Assumption \autoref{ass: boundedness} implies moreover that that for all $n$, $\PP$-almost any $\omega\in\Omega$, $\frac{\nu_0 (m_{1,n})}{\vvvert M_{1,n}\vvvert}>0$. Thus, setting \[d(\omega)=\inf_{n\in \NN} \frac{\nu_0 (m_{1,n})}{\vvvert M_{1,n}\vvvert},\]
$(\nu, d)$ satisfies \eqref{eq : controle d} and $d(\omega)>0$, $\PP(d\omega)$-a.s.
\end{proof}

This provides nice sufficient conditions for \autoref{ass : moments_cd} {\color{black} or \autoref{ass: moments_cd reinforced}} to hold.
\begin{proposition} \label{prop: int logd hennion}
Consider a random, stationary sequence of $p\times p$ matrices $M_n= \left( \mat{M_n}{x,y}\right)_{x,y\in\XX}$, with nonnegative entries, satisfying \autoref{ass: ergodicity}, \autoref{ass: moments_m}, \autoref{ass: hennion}. We assume that there exists a measurable map
\[\omega\in\Omega \mapsto (\nu_\omega,c(\omega)) \in \mathcal{M}_1(\XX)\times [0,1] \] such that 
$\color{black}\PP[c>0]>0$ and \eqref{eq:prop c} hold. Then there exists a random variable $d$ such that $(\nu,c,d)$ is an admissible triplet with $\PP[\gamma>0]>0$. Thus assumption \autoref{ass : moments_cd} holds and Theorem \ref{thme: approx} applies.
\\ If moreover
\begin{enumerate}[i)]
    \item $ \EE\left[-\log c\right]<\infty,$
    \item there exists an deterministic integer $N\in \NN_0$ such that \begin{equation} \label{eq: expected bounded quotients} \EE \left[ \log \sup_{x,y,z\in\XX}\frac{\mat{M_{0,N}}{x,z}}{\mat{M_{0,N}}{y,z}}\right]<\infty, \end{equation}
    \item $\int \left[\log \nu_{\omega}(m_{1,2})\right]^-d\PP(\omega)<\infty$,
\end{enumerate}
Then  \autoref{ass: moments_cd reinforced} also holds.
\end{proposition}
Note that since $\XX$ is finite, $(\delta_j,c)$ satisfies \eqref{eq:prop c} as soon as 
$c\leq \min_{1\leq i \leq p} \frac{M_0(i,j)}{\sum_{l=1}^p M_0(i,l)},$ thus \autoref{ass : moments_cd} holds if with positive probability, $M_0$ has a column with only nonzero coefficients.
\begin{proof}
{\color{black}   By Lemma \ref{lem: pos d hennion}, setting \[d(\omega)=\inf_{n\in \NN} \frac{\nu_\omega(m_{1,n})}{\vvvert M_{1,n}\vvvert},\] $\omega\mapsto (\nu_\omega, d(\omega))$ satisfies \eqref{eq: prop d} and $\PP[d>0]=1$. Noting $\gamma(\omega)=c(\omega)d(\omega)=c_0(\omega)d_1(\omega)$, we notice that $\PP[d_1>0]=\PP[d>0]=1$ and $\PP[c>0]>0$. It holds thus $\PP[\gamma>0]>0$ : \autoref{ass : moments_cd} is satisfied. This proves the first part of the proposition. Let us suppose now that \textit{i)-iii)} hold.} Since $ \EE\left[-\log c\right]<\infty,$ then $\PP[c>0]=1$ and $\PP[\gamma>0]=1$. Let us now prove now that  $\EE |\log\gamma|<\infty$. By the inequality
\[ \vert \log \gamma\vert=-\log(\gamma)\leq-\log(c)-\log(d) \]
and hypothesis
$\displaystyle  \EE[|\log(c)|]<\infty$,  it remains to check that 
\(\EE[|\log(d)|]<\infty.\)
Inequality \eqref{eq : majoration hennion 1} implies that $\PP(d\omega)$-almost surely, for any  $k\geq 1$,
\begin{align*}
-\log d= \log \sup_{n\in\NN} \frac{\vvvert M_{1,n}\vvvert}{\nu_\omega(m_{1,n})} 
&
\leq \max \left( \log \sup_{x,y,z\in \XX} \frac{\mat{M_{1,k+1}}{x,z}}{\mat{M_{1,k+1}}{y,z}}, \max_{1\leq n\leq k} \log \frac{\vvvert M_{1,n}\vvvert}{\nu(m_{1,n})}\right)
\\
&
\leq \log \sup_{x,y,z\in \XX} \frac{\mat{M_{1,k+1}}{x,z}}{\mat{M_{1,k+1}}{y,z}} +\sum_{1\leq n\leq k} \log \frac{\vvvert M_{1,n}\vvvert}{\nu_\omega(m_{1,n})}.
\end{align*}
In particular, setting $k=N$  and applying  condition (\ref{eq: expected bounded quotients}), it holds by stationarity
\[ \EE\left[ \log \sup_{x,y,z\in\XX}\frac{\mat{M_{1,N+1}}{x,z}}{\mat{M_{1,N+1}}{y,z}}\right]<\infty. \]
Consequently, it suffices to prove that  $\displaystyle \EE \left[\log \frac{\vvvert M_{1,n}\vvvert}{\nu_0(m_{1,n})}\right] <\infty$ for any $n \geq 1$. Let us decompose this quantity as
\begin{equation}\label{hzgdv}
 \EE \left[\log \frac{\vvvert M_{1,n}\vvvert}{\nu_0(m_{1,n})}\right] \leq 
\underbrace{\EE \left( \log \vvvert M_{1,n}\vvvert \right)^+ }_{A(n)}
+ 
\underbrace{\EE \left[ \left( \log  \nu_0(m_{1,n}) \right)^- \right]}_{B(n)}. 
\end{equation}

$\bullet$ On the one hand,  the inequality 
$\displaystyle \vvvert M_{1,n}\vvvert \leq \prod_{i=1}^{n-1} \vvvert M_{i,i+1}\vvvert$ readily yields
\[ A(n)\leq \sum_{i=1}^{n-1} \EE\left[ \log(\vvvert M_{i,i+1}\vvvert)^+ \right]=(n-1)\EE\left[ \log(\vvvert M_{0,1} \vvvert)^+\right]<\infty.\]

$\bullet$ On the other hand, for any $x\in\XX$ and $\mathbb P(d\omega)$-a.s. 
\[ m_{1,n}(x)=\delta_x M_{1,2} m_{2,n} \geq c_1 m_{1,2}(x) \nu_{1}(m_{2,n}),\]
where $c_k=c(\theta^k(\omega))$ and $\nu_k=\nu_{\theta^k(\omega)}$.
Consequently, integrating with respect to $\nu_0(dx)$, we obtain
\[  \nu_0(m_{1,n})=\nu_0 M_{1,2} m_{2,n} \geq c_1(\omega)\nu_0(m_{1,2}) \nu_{\theta(\omega)}(m_{2,n})\quad \mathbb P(d\omega){\rm -a.s},\]
which yields, by induction
\[\nu_0(m_{1,n})\geq \prod_{k=1}^{n-1} c_k   \nu_{k-1} (m_{k,k+1}).\]
Consequently, $\mathbb P(d\omega)$-a.s.,  
\[
\left( \log  \nu_0(m_{1,n}) \right)^- \leq  \sum_{k=1}^{n-1}  -\log c_k + \left[\log\nu_{k-1} (m_{k,k+1})\right]^-. 
\]
By stationarity, we deduce
\begin{align*}B(n) &\leq  \sum_{k=1}^{n-1} \EE\left[ -\log c_k \right]+ \EE\left[\left(\log\nu_{k-1} (m_{k,k+1})\right)^- \right]
\\
&= 
(n-1)\EE\left[ -\log c_0 \right]+ \EE\left[\left(\log\nu_{0} (m_{1,2})\right)^- \right]
<\infty.
\end{align*}
Finally, combining these estimates, we get, for any $n\in \NN_0$,
\begin{align*} \EE\left[\log \frac{\vvvert M_{1,n}\vvvert}{\nu_0(m_{1,n})}\right]\leq (n-1)\EE \left[ -\log c_0+ \log\left[\nu_{0}(m_{1,2})\right]^- + \left(\log \Vert m_{0,1}\Vert_\infty \right)^+\right]<\infty.
\end{align*}
\end{proof}

\subsection{Extension in infinite dimension}\label{subs: inFiDi Hilbert}
When $\XX$ is infinite, we need to strengthen the notion of positive matrices as follows.
\begin{definition}
A positive linear map $M$ on $\linf$ is uniformly positive if there exists $K\in \RR_+^*, h\in \mathcal{B}_+(\XX)$, such that, for any $f\in \mathcal{B}_+(\XX)$ there exists $b(f) \in \RR_+$, satisfying 
\[ \frac{1}{K} b(f) h \leq M(f) \leq Kb(f) h .\]
\end{definition}
Notice that when $\XX$ is finite, a matrix of $\accentset{\circ}{\mathcal{S}}$ is uniformly positive. Moreover, in Hennion's work, assumption \autoref{ass: hennion} is used as a sufficient condition to obtain projective contraction properties on the product $M_{k,n}$, with respect to a projective distance called the Hilbert distance (once again, see \cite{busemann_projective_1953,birkhoff_extensions_1957,ligonniere_contraction_2023} for a complement on this distance). In an infinite dimensional setting, this distance can still be defined, and the projective action associated with a positive operator is contracting if and only if the operator is uniformly positive (a proof of this claim is proposed in \cite{ligonniere_contraction_2023}). Uniform positivity is therefore the appropriate infinite dimensional generalization of positivity in our context, and condition \autoref{ass: hennion} can thus naturally be replaced with the restrictive condition 
\begin{manualass}{\begin{NoHyper}\ref{ass: hennion}\end{NoHyper}'} \label{ass: uniform_pos}
    \(\PP\left( \bigcup_{n\in \NN} \left \{ M_{0,n} \gg 0 \right\} \right)=1,\)
\end{manualass}
where we note $M\gg 0$ iff $M$ is uniformly positive.

The present subsection aims at comparing our result with the natural extensions of Hennion's work to infinite dimensional settings. For this purpose, the following Lemma extends the idea of Lemma \ref{lem: pos d hennion} to an infinite dimensional setup, assuming \autoref{ass: uniform_pos} instead of \autoref{ass: hennion}.
\begin{lemma} \label{lem: unif pos X infini}
Consider a random stationary sequence of elements of $\MX$, satisfying \autoref{ass: boundedness} and \autoref{ass: uniform_pos}. Then for any measurable map $\omega\in\Omega\mapsto \nu(\omega)\in \mathcal{M}_1(\XX)$, there exists a random variable $d$ such that $\PP[d>0]=1$ and \eqref{eq: prop d} holds.
\end{lemma}
\begin{proof}
For almost any $\omega$ and any $1\leq k\leq n$, $m_{k,n}\in \linf$, it holds,
\[ \frac{m_{1,n}(x)}{\nu(m_{1,n})}=\frac{\delta_x M_{1,k}{m_{k,n}}}{\nu_0 M_{1,k} m_{k,n}} \leq \sup_{y\in \XX, f\in \linf } \frac{ \delta_y M_{1,k} f }{\nu_0 M_{1,k} f}.\]
Taking a supremum in $x\in\XX$, we get, for any $k\leq n$,
\beq \label{eq : csq unifpos 1} \frac{\vvvert M_{1,n}\vvvert}{\nu_0(m_{1,n})} \leq \sup_{y\in \XX,f\in \B(\XX)} \frac{ \delta_y M_{1,k} f }{\nu_0 M_{1,k} f}.\eeq
By \autoref{ass: uniform_pos}, let $k_\omega$ be a random integer  such that $\PP(d\omega)$-almost surely, $M_{1,k_\omega}(\omega)\gg 0$. Then, almost surely, there is $K \in \RR_+^*$, $h\in \mathcal{B}_+(\XX) $ such that for any $f\in \mathcal{B}_+(\XX)$, there exists $b(f)\geq0$, satisfying
\beq \label{eq: unifpos} K^{-1} b(f) h \leq M_{1,k_\omega}f \leq K b(f) h.\eeq
From \eqref{eq: unifpos}, we deduce $K^{-1} m_{1,k_\omega}(x)\leq h(x)b(\mathds{1}) \leq K m_{1,k_\omega}(x)$. By \autoref{ass: boundedness}, $m_{1,k_\omega}$ is a bounded and positive function, thus so is $h$. Moreover, $b(\mathds{1})>0$, $\nu_0(m_{1,k})\leq K\nu_0(h)b(\mathds{1})$, thus $\nu(h)>0$.
Therefore, for any $x\in\XX$, any $f\in\mathcal{B}_+(\XX)$, it holds $\PP(d\omega)$ almost surely:
\begin{equation}\label{eq: csq unifpos 2} \frac{\delta_x M_{1,k_\omega}f }{\nu_0 M_{1,k_\omega} f} \leq K^2 \frac{h(x)}{\nu_0(h)}\leq K^3\frac{\vvvert M_{1,k_\omega} \vvvert}{b (\mathds{1}) \nu_0(h)}\leq K^4 \frac{\vvvert M_{1,k_\omega} \vvvert}{\nu_0(m_{1,k_\omega})}.\end{equation}
Finally, combining \eqref{eq: csq unifpos 2} with \eqref{eq : csq unifpos 1}, we get for almost any $\omega$,
\[\limsup_n \frac{\vvvert M_{1,n} \vvvert}{\nu_0(m_{1,n})}\leq \sup_{n\geq k_\omega} \frac{\vvvert M_{1,n} \vvvert}{\nu_0(m_{1,n})} \leq K^4 \frac{\vvvert M_{1,k_\omega} \vvvert}{\nu(m_{1,k_\omega})} <\infty.\]
Since moreover almost surely, for all $n\geq 1$, $\frac{\vvvert M_{1,n} \vvvert}{\nu_0(m_{1,n})}$ is finite, then almost surely,
\[\sup_{n\in\NN} \frac{\vvvert M_{1,n} \vvvert}{\nu_0(m_{1,n})} <\infty.\]
Let us set 
\[d(\omega)=\inf_{n\geq 1} \frac{\nu_0(m_{1,n})}{\vvvert M_{1,n} \vvvert }.\]
Then $d$ clearly satisfies \eqref{eq: prop d} and $\PP[d>0]=1$.
\end{proof}
The uniform positivity property is interesting to deal with many systems where the mass is sufficiently well mixed. We illustrate this on the following example.
\begin{example}
Take $\XX=[0,1]^d$, choose an ergodic dymanical system $(\Omega, \mathcal{A},\PP,\theta)$ and associate with each $\omega\in\Omega$ a bounded measurable function $m_\omega:\XX\rightarrow \RR_+$, and a continuous function $Q_\omega:\XX^2\rightarrow \RR_+^*$, such that $\int_\XX Q_\omega(x,y)dy=1$ for all $x$ and all $\omega$.
We define $M$ by setting, for each $\omega\in\Omega$ and $f\in \B(\XX),$
$$M(\omega)(f):x\mapsto m_\omega(x)\int_{\XX} f(y)Q_\omega(x,y)dy.$$
Notice that $m_\omega(x)=\Vert \delta_x M(\omega)\Vert$.
Our system clearly satisfies \autoref{ass: ergodicity} and \autoref{ass: boundedness}. Assumptions \autoref{ass: moments_m} and \autoref{ass: moments_m reinforced} just translate into a log-integrability assumption on $\omega\mapsto \Vert m_\omega \Vert_\infty=\vvvert M(\omega)\vvvert $.
\\ In terms of populations, this model can represent the spatial evolution of a population in the compact domain $\XX$. The quantity $m_\omega(x)$ represents the size of the offspring of an individual located at $x$, and the kernel $Q_\omega(x,y)$ represents the dispersion of its offspring in the domain $\XX$. The dependence in $\omega$ of these quantities models an time inhomogeneity of these quantities.
$M(\omega)$ is clearly uniformly positive, for each $\omega$, since
\begin{equation*}K^{-1}m_\omega \mu(f) \leq M(\omega)(f)\leq K m_\omega \mu(f),\end{equation*}
where $\mu$ refers to the Lebesgue measure on $\XX$, and $$K=K_\omega=\max\left(\left(\inf_{u,v} Q_\omega(u,v)\right)^{-1}, \sup_{u,v} Q_\omega(u,v)\right)<\infty$$ since $\XX$ is compact and $Q_\omega$ is continuous and positive.
Setting $\nu_\omega=\mu$ and $c_\omega=K_\omega^{-1}$ for all $\omega$, the left hand side inequality implies \eqref{eq:prop c}. Moreover, since $M(\omega)$ is uniformly positive, \autoref{ass: uniform_pos} holds, thus by Lemma \ref{lem: unif pos X infini} there exists a random variable $d$ such that $\PP[d>0]>0$ and $\omega\mapsto(\nu_\omega, d(\omega))$ satisfy \eqref{eq: prop d}. Therefore \autoref{ass : moments_cd} holds and Theorem \ref{thme: approx} applies.
\\ One can notice additionally that in this context, each matrix of the product is almost surely uniformly positive, therefore the proof of Lemma \ref{lem: unif pos X infini} yields the explicit control
\[d(\omega)> \frac{\int M_{1,2}(\mathds{1})(x)dx}{K_\omega^4\sup_{x\in \XX}M_{1,2}(\mathds{1})(x)} =\frac{\int m_{\theta(\omega)} (x)dx}{K_\omega^4\Vert m_{\theta(\omega)} \Vert_\infty}.\]
Thus \autoref{ass: moments_cd reinforced} reduces to a $\log$-integrability condition both on the coefficient $K_\omega$ and on the quotient $\frac{\int m_\omega (x)dx}{\Vert m_\omega \Vert_\infty}$.
\end{example}
\smallskip This example illustrates how Proposition \ref{prop: int logd hennion} from the previous subsection can be adapted, replacing \autoref{ass: hennion} by \autoref{ass: uniform_pos}. To tackle the integrability of $\log \gamma$, one can replace \eqref{eq: expected bounded quotients} by 
\begin{equation} \label{eq : expected bounded quotients infdim} \EE\left[ \log \sup_{x\in\XX, f\in \linf}\frac{\delta_x M_{1,N}f }{\nu_0 M_{1,N} f} \right]<\infty. \end{equation}
This yields a counterpart of Proposition \ref{prop: int logd hennion} in a infinite dimensional setup.

%% file: 5_branching.tex
\section{Application to products of infinite Leslie matrices}\label{sec: branching models}
The previous section focused on products of matrices with positive entries, and more generally, products of uniformly positive operators. This kind of products can be efficiently studied with methods based on projective contractions relatively to the Hilbert metric. The main interest of our techniques, based on Doeblin contractions, is their application to products of operators which are not uniformly positive. The goal of this section is to illustrate how such products can be studied with our theorems. We have chosen to focus here on a quite simple but natural example with no uniform positivity properties : the infinite Leslie matrices.

\subsection{Introduction to Leslie matrices}
In this section, we set $\XX=\NN_0$, thus the operators of $\MX$ can be represented as infinite matrices. We choose to consider infinite Leslie matrices, which have the following form : for any $\omega\in \Omega$,
 \begin{equation} M(\omega)=\begin{pmatrix} \label{eq : Leslie_PF}
        f_0(\omega) & s_0(\omega) &  0  &  0  & \dots  \\
        f_1(\omega) &  0  & s_1(\omega) &  0  & \dots  \\
        f_2(\omega) &  0  &  0  & s_2(\omega) & \ddots \\
        f_3(\omega) &  0  &  0  &  0  & \ddots \\
        \vdots & \vdots & \vdots & \vdots & \ddots \\
 \end{pmatrix}. \end{equation}
 where the entries $\left(f_k(\omega)\right)_{k\in\NN_0}$, $\left(s_k(\omega)\right)_{k\in\NN_0}$ are nonnegative and $\sup_{x\in \XX} s_x(\omega)+f_x(\omega)<\infty$. Notice that such a matrix is not uniformly positive, since there are zeros on every row and every column but the first one. Moreover, if $Q$ is a product of $k$ matrices of this shape, the $(x,y)$-entry $\left[Q\right]_{x,y}=0$ whenever $y\geq x+k+1$. This prevents any product of such matrices from being uniformly positive. This example is therefore a typical situation where \autoref{ass: uniform_pos} does not hold. 

Such matrices appear naturally when studying the dynamics of a population counting individuals according to their age. The coefficients $f_x$ (respectively $s_x$) represent the mean number of individuals of age $0$ (respectively of age $x+1$) created by an individual of age $x$, that is the mean size of the offspring of an individual of age $x$ (respectively the survival rate of individuals of age $x$). Usually, only a finite number of age classes are defined, thus $\XX= \llbracket 0,p \rrbracket$, and one considers finite versions of such matrices, called Leslie matrices, see for example \cite{caswell_life_2010}. However, it is natural to extend their definition to an infinite number of age classes ($\XX=\NN_0$) obtaining infinite matrices with this shape. Indeed, several articles already study age-structured populations with an unbounded set of possible ages, see e.g. \cite{bansaye_ergodic_2020, jasinska_dynamics_2022, oelschlager_limit_1990}. Therefore, products of random matrices shaped as in \eqref{eq : Leslie_PF} model the dynamics of an age structured population evolving in a randomly changing environment which affect their reproductive behavior. This is the kind of matrices we are studying in this section.
Let us note from now on \[s_x^k(\omega)=s_x\circ\theta^k(\omega)\text{ and } f_x^k(\omega)=f_x\circ\theta^k(\omega),\] so that $(s_x^k,f_x^k)_{x\in\XX}$ are the nonzero entries of the random matrix $M_k(\omega)=M\circ\theta^k$.
We introduce the quantities
\[d'(\omega)=\sup_{k\in\NN_0, x\leq y\in \XX} \frac{f_{y}^k(\omega)}{f^k_x( \omega)}\geq 1,\]
and
\[d''(\omega)= \sup_{x\in \XX, k\in\NN_0}\frac{s_x^0(\omega)\cdots s_{x+k}^{k}(\omega)}{s_0^0(\omega)\dots s_k^k(\omega)}\geq 1,\]
which are useful to construct an admissible triplet.
\subsection{Ergodic behavior of products of random Leslie matrices} \label{subs: assumptions_age}
The following proposition provides sufficient conditions for assumptions \autoref{ass: moments_m} and \autoref{ass : moments_cd} to hold in the case of products of infinite Leslie matrices.
\begin{proposition}\label{prop:hyp_leslie}
Consider a random matrix product with $\XX=\NN_0$ and suppose that for any $\omega\in \Omega$, $M(\omega)$ is of the form of equation \eqref{eq : Leslie_PF}, with $\sup_{x\in\XX} s_x(\omega)+f_x(\omega) <\infty$. Suppose that \autoref{ass: ergodicity} is satisfied, and $\PP(d\omega)$-almost surely, it holds
\begin{enumerate}
    \item[i)] $f_x(\omega)+s_x(\omega)>0$ for all $x\in\XX$ ;
    \item[ii)] $\EE\left[\log^+ \left(\sup_{x\in\XX} s_x+f_x\right) \right]<\infty$ ;
\end{enumerate}
then \autoref{ass: boundedness} and \autoref{ass: moments_m} hold. {\color{black} If moreover
\begin{enumerate}
    \item[iii)] There exists a deterministic real $A>0$ such that $\PP(d\omega)$-almost surely, $\sup_{x\leq y } \frac{f_{y}^0}{f_x^0}\leq A$,
    \item[iv)] $ \PP[ \sup_{x\in\XX}\frac{s^0_x}{f^0_x}<\infty, d^{''}\circ \theta <\infty]>0$
\end{enumerate}
then Assumption \autoref{ass : moments_cd} holds, and so do the conclusions of Theorem \ref{thme: approx}.
\\ Finally, if additionally }
\begin{enumerate}
    \item[v)] $\PP(d\omega)$-almost surely, $\sup_{x\in\XX}\frac{s_x}{f_x}<\infty $ and $\EE\left[ \log^+(\sup_{x\in\XX}\frac{s_x}{f_x})\right]<\infty $,
    \item[vi)] $\PP(d\omega)$-almost surely, $d''(\omega)<\infty$ and $\EE \left|\log d''\right|<\infty,$
\end{enumerate}
then $M$ satisfies also Assumption \autoref{ass: moments_cd reinforced}. 
\end{proposition}

In Proposition \ref{prop:hyp_leslie}, we've reduced Assumptions \autoref{ass: boundedness} to \autoref{ass : moments_cd} to a series of conditions on the law of the coefficients of the random matrix $M_0$, together with finiteness and integrability conditions on $d''$. The hardest conditions to check are the ones involving $d''$, since checking them requires to consider the joint law of all the $M_{0,n}$ and not only $M_0$. We were not able to find a general sufficient condition for the positivity and $\log$-integrability of $d''$. However, we provide the following quite restrictive sufficient condition.
\begin{remark}\label{rem:hyp leslie}
Consider a random, stationary sequence of matrices of the form of equation \eqref{eq : Leslie_PF}, and assume that there exists an integer $x_0\in\XX$, such that almost surely, the sequence $(s_x(\omega))_{x\geq x_0}$ is non increasing. Suppose also that almost surely, for all $x\leq x_0$, $s_x>0$.
Then, $\PP(d\omega)$-almost surely
\[ d''\leq \left(\sup_{i\leq x_0} \sup_{x \leq y\leq x_0} \frac{s_y^i}{s_x^i}\right)^{x_0}<\infty.\]
Moreover, if $\EE \left|\log \frac{s_y^i}{s_x^i}\right|<\infty$ for any $x\leq y\leq x_0$, then $\EE |\log d''|<\infty$.
\end{remark}
In the context of an age structured population, $s_x$ represents the frequency of individuals of age $x$ surviving to the next time step, and thus being replaced by individuals of age $x+1$. Assuming that $(s_x(\omega))_{x\geq x_0}$ is decreasing implies that the older individuals get, the more they tend to die, which is a reasonable assumption. However this condition is somewhat unsatisfying in a more general setting.

We split the proof of Proposition \ref{prop:hyp_leslie} into several lemmas that involve different groups of assumptions.
Notice first that most quantities involved in Assumptions \autoref{ass: boundedness} to \autoref{ass : moments_cd} are explicit in terms of the $(f_x,s_x)$. Indeed :
\begin{lemma} \label{lem:explicitquantites}
Consider a product of stationary random Leslie matrices, in the form of equation \eqref{eq : Leslie_PF}. Then \autoref{ass: boundedness} and \autoref{ass: moments_m}  are satisfied if and only if all the following conditions hold simultaneously :
\begin{itemize}
    \item $\PP(d\omega)$-almost surely, for each $x\in\XX$, $f_x(\omega)>0$ or $s_x(\omega)>0$
    \item $\EE\left[\log^+ \left(\sup_{x\in\XX} f_x+s_x\right) \right]<\infty$.
\end{itemize}
\end{lemma}
\begin{proof}
This lemma is straightforward after noticing that for any $x\in\XX,\omega\in\Omega$, 
\[m_{0,1}(x,\omega)=f_x(\omega)+s_x(\omega).\]
\end{proof}
Moreover, in this model, \eqref{eq:prop c} is well behaved and it is quite clear how to construct a non trivial couple $(\nu,c)$.
\begin{lemma} \label{lem:controle c}
Consider a product of stationary, random Leslie matrices and assume that \autoref{ass: boundedness} holds. Consider a map $\omega\mapsto (\nu_\omega,c(\omega))$. Then for almost any $\omega$ such that $\nu_\omega\neq \delta_0$, \eqref{eq:prop c} implies
\( c(\omega)=0.\)
If $\nu_\omega=\delta_0$, then \eqref{eq:prop c} is equivalent to
\[c(\omega)\leq \inf_{x\in\XX} \frac{f_x(\omega)}{f_x(\omega)+s_x(\omega)}= \left(1+\sup_{x\in\XX}\frac{s_x(\omega)}{f_x(\omega)}\right)^{-1}.\]

\end{lemma}
\begin{proof}
Notice that for any $x\in\XX$, $\omega\in\Omega$,
\[ \delta_x M_{0,1} =f_x(\omega)\delta_0+ s_x(\omega)\delta_{x+1}.\]
Let $\omega\in\Omega$ such that $\nu_\omega\neq \delta_0$. Then, there exists $k>0$ such that $\nu_\omega(k)>0$. In particular \eqref{eq:prop c} implies
\[ 0=f_k(\omega)\delta_0(\mathds{1}_k)+ s_k(\omega)\delta_{k+1}(\mathds{1}_k)\geq c(\omega) m_{0,1}(k)\nu(k),\]
By \autoref{ass: boundedness}, almost surely, $m_{0,1}(k)>0$, which implies that $c(\omega)=0$.
Conversely, if $\nu_\omega=\delta_0$, \eqref{eq:prop c} implies
\[f_x(\omega)\delta_0+ s_x(\omega)\delta_{x+1}\geq c(\omega) m_{0,1}(x)\delta_0,\]
which is equivalent to 
\[f_x(\omega)\geq c(\omega) m_{0,1}(x)=c(\omega)(f_x(\omega)+s_x(\omega))\]
for all $x\in \XX$. This yields the desired result. 
\end{proof}
As a consequence, we set from now on $\nu_\omega=\delta_0$ and $c_\omega=\left(1+\sup_{x\in\XX}\frac{s_x(\omega)}{f_x(\omega)}^{-1}\right)^{-1}$. Assumption \textit{iv)} of Proposition \ref{prop:hyp_leslie} guarantees that $\PP[c>0]>0$.
Let us try now to exhibit a random variable $d$ such that $(\nu=\delta_0,d)$ satisfy \eqref{eq: prop d}.

\begin{lemma}\label{lem:controle d}
Consider a product of stationary random Leslie matrices, of the form of equation \eqref{eq : Leslie_PF}. Set $$d(\omega)=\frac{1}{d'\circ\theta(\omega)d''\circ\theta(\omega)}.$$
Then $(\delta_0,d)$ satisfy \eqref{eq: prop d}.
\end{lemma}

\begin{proof}
Let $n \geq 1, x\in \XX$, $\omega\in\Omega$, it holds 
\[m_{1,n}(x)=\delta_x M_0\cdots M_{n-1} \mathds{1}=\sum_{i_0,i_1,\cdots i_n\in \NN_0} \delta_x(i_0)\mat{M_1}{i_0,i_1} \cdots \mat{M_{n-1}}{i_{n-1},i_n}.\]
Thus
\[m_{1,n}(x)=\sum_{i_1,\cdots i_n\in \NN_0} \mat{M_1}{x,i_1} \cdots \mat{M_{n-1}}{i_{n-1},i_n}.\]
Let us rearrange this sum according to the first index $k\leq n$ such that $i_k=0$ :
\[ \begin{split} m_{1,n}(x)=& \sum_{k=1}^n \sum_{i_1,\cdots i_{k-1}>0}\mat{M_0}{x,i_1} \cdots \mat{M_{k-1}}{i_{k-1},0}  \sum_{i_{k+1},\cdots i_n\in \NN_0} \mat{M_{k}}{0,i_{k+1}}\cdots \mat{M_{n-1}}{i_{n-1},i_n}
\\ +& \sum_{i_1,\cdots i_n>0} \mat{M_0}{x,i_1} \cdots \mat{M_{n-1}}{i_{n-1},i_n}.\end{split}\]
Notice that 
\[\sum_{i_{k+1},\cdots i_n\in \NN_0} \mat{M_{k}}{0,i_{k+1}}\cdots \mat{M_{n-1}}{i_{n-1},i_n}=m_{k,n}(0).\]
Moreover, the matrices $M_k$ are shaped according to \eqref{eq : Leslie_PF}. Thus for any $i\geq 0,j>0$, in order for $\mat{M_{k}}{i,j}$ to be non zero, one must have $j=i+1$. Thus :
\begin{align*}m_{0,n}(x)=\sum_{k=1}^n &\mat{M_0}{x,x+1} \cdots \mat{M_{k-1}}{x+k-1,0}  m_{k,n}(0)\\ &+\mat{M_0}{x,x+1} \cdots \mat{M_{n-1}}{x+n-1,x+n}.\end{align*}
Therefore
\[m_{0,n}(x)=\sum_{i=0}^{n-1}s_x^0 s_{x+1}^{1}\dots s_{x+i-1}^{i-1} f_{x+i}^{i} m_{i+1,n}(0) + s_x^0\dots s_{x+n-1}^{n-1}.\]
This is true in particular for $x=0$ :
\[m_{0,n}(0)=\sum_{i=0}^{n-1}s_0^0 s_{1}^{1}\dots s_{i-1}^{i-1} f_{i}^{i} m_{i+1,n}(0) + s_0^0\dots s_{n-1}^{n-1}.\]
By definition of $d',d''$, it holds, for any $k \in \NN_0$ and any $x\in\NN_0$,
\[f_{x+i}^{i}\leq d'f_{i}^{i},\]
and\[s_x^0 s_{x+1}^{1}\dots s_{x+i-1}^{i-1}\leq d''s_0^0 s_{1}^{1}\dots s_{i-1}^{i-1}.\]
Therefore, controlling independently each term of the sum yields
\[m_{0,n}(x)\leq d'd'' \sum_{i=0}^{n-1}s_0^0 s_{1}^{1}\dots s_{i-1}^{i-1} f_{i}^{i} m_{i+1,n}(0)+ d'' s_0^0\dots s_{n-1}^{n-1}\leq d'd'' m_{0,n}(0),\]
thus for all $n\geq 0$
\[\frac{1}{d'(\omega)d''(\omega)}\leq \inf_{x\in\XX} \frac{m_{0,n}(0)}{m_{0,n}(x)}.\]
By stationarity,
\[\frac{1}{d'(\theta(\omega))d''(\theta(\omega))}\leq\inf_{n\geq 1}\frac{m_{1,n}(0)}{\Vert m_{1,n} \Vert_\infty} .\]
As a consequence, setting $d(\omega)=\left(d'(\theta(\omega))d''(\theta(\omega))\right)^{-1}$  is enough for $\omega\mapsto (\delta_0,d(\omega))$ to satisfy \eqref{eq: prop d}.
\end{proof}
Let us focus on $d'(\omega)$.
\begin{lemma}\label{lem:controle fertilites}
Consider a random product of matrices of the form of equation \eqref{eq : Leslie_PF}, satisfying \autoref{ass: ergodicity}. Then the random variable $d'$ is $\PP(d\omega)$-almost surely finite if and only if there exists $A>0$ such that 
\begin{equation}\label{eq: bounded growth fertilites} \PP\left[ \sup_{x\leq y } \frac{f_{y}^0}{f_x^0} \leq A \right]=1.\end{equation} In this case $d'(\omega)\leq A$, $\PP(d\omega)$-almost surely.
{\color{black} If \eqref{eq: bounded growth fertilites} fails, then $d'(\omega)=+\infty$, $\PP(d\omega)$- almost surely.}
\end{lemma}
\begin{proof}
Notice that \[d'(\omega)=\sup_{k\in\NN_0} \sup_{x\leq y } \frac{f_{y}^k}{f_x^k}=\sup_{k\in\NN_0} X\circ \theta^k(\omega).\]
where \[X(\omega)=\sup_{x\leq y } \frac{f_{y}(\omega)}{f_x(\omega)}=\sup_{x\leq y } \frac{f_{y}^0}{f_x^0}.\]
Since $\theta$ is an ergodic mapping, $\sup_{k\in\NN_0} X\circ\theta^k$ is $\PP(d\omega)$-almost surely equal to the supremum of the support of $X$. In particular $\sup_{k\in\NN_0} X\circ\theta^k$ is finite almost surely if the support of $X$ is bounded. {\color{black} Conversely, if the support of $X$ is unbounded, then $\sup_{k\in\NN_0} X\circ\theta^k$ is infinite almost surely.}
\end{proof}

Putting these lemmas together allows to prove Proposition \ref{prop:hyp_leslie}.
\begin{proof}[Proof of Proposition \ref{prop:hyp_leslie}]
The assumptions \textit{i), ii)} of Proposition \ref{prop:hyp_leslie} are the conditions mentioned in Lemma \ref{lem:explicitquantites}. Hence, this lemma allows to check \autoref{ass: boundedness} and \autoref{ass: moments_m}.
We set now for any $\omega$, $\nu_\omega=\delta_0$
\[c(\omega)=\left( 1+ \sup_{x\in\XX} \frac{s_x(\omega)}{f_x(\omega)}\right)^{-1}\]
and 
\[d(\omega)=\frac{1}{d'(\theta(\omega))d''(\theta(\omega))}.\]
{\color{black} Lemma \ref{lem:controle c} and  Lemma \ref{lem:controle d} respectively guarantee that \eqref{eq:prop c} and \eqref{eq: prop d} are satisfied.
By Lemma \ref{lem:controle fertilites}, assumption \textit{iii)} guarantees that $d'(\omega)<\infty, \PP(d\omega)-$almost surely, and by stationarity, the same holds for $d'\circ \theta(\omega)$.
As a consequence, since by \textit{iv)} $\PP[d''\circ \theta>0,c_0>0]>0$, it holds with positive probability $\gamma(\omega)=c(\omega)d(\theta(\omega))=c(\omega)\left(d'(\theta(\omega))d''(\theta(\omega))\right)^{-1}>0$ thus $\PP[\gamma>0]>0$. This proves that \textit{i)-iv)} imply \autoref{ass : moments_cd}.
\\ If moreover assumptions \textit{v), vi)} hold, then $d''\circ \theta$ and $\sup_{x\in\XX} \frac{s_x(\omega)}{f_x(\omega)}$ are finite $\PP(d\omega)-$a.s., thus $\gamma>0$ almost surely.}
In this case it also holds
\[\EE|\log\gamma|\leq \EE\left[-\log c\right] +\EE\left[-\log d\right].\]
On the one hand, 
\[\EE \left[-\log c\right]=\int \log\left(1+ \sup_{x\in\XX} \frac{s_x(\omega)}{f_x(\omega)}\right)d\PP(\omega).\]
Notice that for any positive real variable $X$, $\log(1+X)$ is integrable as soon as $\log(X)^+$ is integrable. Since we've assumed that \[\int \log\left(\sup_{x\in\XX} \frac{s_x(\omega)}{f_x(\omega)}\right)^+d\PP(\omega)<\infty,\]
then \[\EE[-\log c]<\infty.\]
On the other hand,
\[\EE\left[-\log d\right] \leq \EE \log d' + \EE \log d''.\]
Since $\sup_{x\leq y}\frac{f_{y}^0}{f_x^0}\leq A$ almost surely, then by stationarity, almost surely, 
\[ 1\leq d'= \sup_{k\in\NN_0} \sup_{y\geq x}\frac{f_{y}^k}{f_x^k} \leq A.\]
Thus $\log d'$ is bounded and integrable. We have assumed additionally that $\log d''$ was integrable. This is enough to conclude to the integrability of $|\log \gamma|$, which proves assumption \autoref{ass: moments_cd reinforced}.
\end{proof}

\subsection{A situation where \texorpdfstring{$\gamma=0$}{gamma=0}} \label{subs : contrex}
It was not clear to us how strong an assumption \autoref{ass : moments_cd} is, or whether it was hard to find a system breaking it while satisfying all the other assumptions. We shall present here an example of an infinite Leslie matrix, such that $\gamma=0$ even if all other assumptions are satisfied. This example is in a deterministic environment, that is $|\E|=1$, $\Omega=\E^{\NN_0}$, $|\Omega|=1$.
The random matrix $M(\omega)$ is therefore constant, and $M_{0,n}=M^n$. Let us set :
\begin{equation}\label{eq : Leslie 2} M=\begin{pmatrix}
        \delta m(0) & \bc m(0) &  0  &  0  & \dots  \\
        \delta m(1) &  0  & \bc m(1) &  0  & \dots  \\
        \delta m(2) &  0  &  0  & \bc m(2) & \ddots \\
        \delta m(3) &  0  &  0  &  0  & \ddots \\
        \vdots & \vdots & \vdots & \vdots & \ddots \\
    \end{pmatrix}. \end{equation}
    where $\delta\in (0,1)$, and $m(x)=\delta m(x)+\bc m(x)$ is the mean offspring size of an individual of age $x$.
    Such a model satisfies \autoref{ass: boundedness} and \autoref{ass: moments_m}, as soon as $x\mapsto m(x)$ is bounded and positive, since $\delta>0$. The ergodicity and integrability properties are trivially satisfied since this model is in a constant environment. Moreover, Lemma \ref{lem:controle c} applies, therefore setting $\nu=\delta_0$, and $c= \left(1+\sup_{x\in\XX} \frac{s_x}{f_x}\right)^{-1}=\delta$, the couple $(\nu,c)$ satisfies \eqref{eq:prop c}. Let us prove that we can tune the parameters $x\mapsto m(x)$ and $c=\delta$ in such a way that the only $d$ satisfying \eqref{eq: prop d} is $d=0$.  
   
Consider a sequence of integers $(\varepsilon_x)_{x\in\NN_0}\in \{0,1\}^{\NN_0}$, such that :
    \begin{itemize}
        \item There are arbitrarily long subsequences of  consecutive $1$ in the sequence $(\varepsilon_x)$.
        \item Noting $S_x=\sum_{k=0}^{x-1} \varepsilon_k$, $\frac{S_x}{x}\longrightarrow 0$ as $x\rightarrow \infty$.
        \item There exists $\alpha<1$ such that for all $x\in \NN_0$, $\frac{S_x}{x}\leq \alpha$.
    \end{itemize}
    Let $a$ be a real number such that $a>1$. Then, we set, for any $x\in\XX$, \[ m(x)= 1+(a-1)\varepsilon_x.\] 
Defined as such, $m$ is a positive and bounded function, thus Assumptions \autoref{ass: boundedness} and \autoref{ass: moments_m} are satisfied.
This yields that for any sequence $(x_i)$, 
\[\prod_{i=0}^{n-1}m(x_i)=a^{\sum_{i=0}^{n-1} \varepsilon_{x_i}}.\]
Moreover, since this model is in constant environment, $M_{0,n}=M_{1,n+1}=M^n$, thus 
\begin{align*}m_{1,{n+1}}(x)&=m_{0,n}(x)=\sum_{\substack{x_0=x,\dots x_n\in\NN \\ x_{i+1}\in \{x_i+1,0\}}} \delta^{N(x_0,\dots x_{n})} \bc^{n-N(x_0,\dots x_{n})}\prod_{i=0}^{n-1} m(x_i)\\ &\geq \bc^n\prod_{i=0}^{n-1}m(x+i),\end{align*}
where $N(x_0,\dots x_n)=\vert \{ 1\leq i\leq n | x_i=0\}\vert.$ Then
\[ m_{0,n}(x) \geq \bc^n a^{\sum_{i=x}^{x+n-1}\varepsilon_i},\]
In particular, $x$ can be chosen such that $\varepsilon_x=\cdots \varepsilon_{x+n-1}=1$, which implies that
\[ \Vert m_{0,n}\Vert_\infty\geq m_{0,n}(x) \geq (a\bc)^{n}.\]
On the other hand
\[ m_{0,n}(0)\leq 2^{n} \sup_{\substack{x_0=0,\dots x_n\in\NN_0 \\ x_{i+1}\in \{x_i+1,0\}}}\prod_{i=0}^{n-1}m(x_i)\leq 2^{n} a^{\sup_{\substack{x_0=0,\dots x_n\in\NN_0 \\ x_{i+1}\in \{x_i+1,0\}}} \sum_{i=0}^{n-1} \varepsilon_{x_i}}.\]
A sequence $(x_i)_{0\leq i \leq n}$ of integers such that  $x_0=0$ and for each $i$, $x_{i+1}\in \{x_i+1,0\}$ is entirely determined by the sequence $(T_k)_{k}$ of the lengths of its excursions away from zero. By convention, if there are only $p$ excursions away from zero, we set $T_{p}$ such that $T_0+\cdots+T_p=n$ and $T_{p+1}=\cdots=T_n=0.$ The $(x_{T_0+\cdots+T_{i-1}})_{i\leq p}$ are the only zero terms in the sequence $(x_1,\dots x_n)$, and $T_0+\cdots + T_n\leq n-1$.
Thus 
\[\sup_{\substack{x_0=0,\dots x_n\in\NN \\ x_{i+1}\in \{x_i+1,0\}}} \sum_{i=0}^{n-1} \varepsilon_{x_i}\leq \sup_{T_0+\dots T_n=n} \sum_{i=0}^{n} S_{T_i}\leq \alpha \sum_{i=0}^n T_i\leq \alpha n,\]
and 
\[m_{0,n}(0)\leq (2a^\alpha)^n.\]
Hence
\[\frac{\Vert m_{1,n+1} \Vert_\infty}{m_ {1,n+1}(0)}=\frac{\Vert m_{0,n} \Vert_\infty}{m_ {0,n}(0)}\geq \left(\frac{a\bc}{2a^\alpha}\right)^n=\left(\frac{a^{1-\alpha}(1-c)}{2}\right)^{n}\underset{n\rightarrow \infty}{\longrightarrow} \infty,\]
whenever $\frac{a^{1-\alpha}(1-c)}{2}>1$. For \eqref{eq: prop d} to hold, we must have $d\leq \inf_{n\geq 1}\left(\frac{m_{1,n}(0)}{\Vert m_{1,n} \Vert_\infty}\right)$. Thus for any values of $\alpha,\delta\in(0,1)$, if $a$ is large enough, then \eqref{eq: prop d} implies $d=0$.